\documentclass{article}
\usepackage[%
journal=DM,    
lang=american,   
]{ems-journal}

\newcommand{\cI}{{\mathcal I}}
\usepackage{hyperref}
\usepackage{tikz}
\usepackage{tikz-cd}
\usetikzlibrary{positioning}

\newtheorem{theorem}{Theorem}[section]
\newtheorem{theoremint}{Theorem}

\newtheorem{Definition}[theorem]{Definition}
\newtheorem{Example}[theorem]{Example}
\newtheorem{Lemma}[theorem]{Lemma}
\newtheorem{Corollary}[theorem]{Corollary}
\newtheorem{Proposition}[theorem]{Proposition}
\newtheorem{Remark}[theorem]{Remark}

\newcommand{\la}{\langle}
\newcommand{\ra}{\rangle}
\newcommand{\exd}{\mathrm{d}}
\newcommand{\R}{{\mathbb R}}
\newcommand{\N}{{\mathbb N}}
\newcommand{\cL}{\mathcal{L}}
\newcommand{\p}{\mathrm{p}}
\newcommand{\pderiv}[2]{\frac{\partial #1}{\partial #2}}

\numberwithin{equation}{section}

\begin{document}

\title{Singular Riemannian foliations \\ and $\mathcal{I}$-Poisson manifolds}
\titlemark{Singular Riemannian Foliations and $\mathcal{I}$-Poisson manifolds}



\emsauthor{1}{
	\givenname{Hadi}
	\surname{Nahari}
	\orcid{0009-0006-4667-6999}}{H.~Nahari}
\emsauthor{2}{
	\givenname{Thomas}
	\surname{Strobl}
	\mrid{334505}
	\orcid{0000-0001-7984-7172}}{T.~Strobl}

\Emsaffil{1}{
	\department{1}{Laboratoire d’Analyse et de Mathématiques Appliquées}
	\organisation{1}{Universit\'e Paris Est Cr\'eteil}
	\rorid{1}{05ggc9x40}
	\address{1}{61 Av. du Général de Gaulle}
	\zip{1}{94000}
	\city{1}{Créteil}
	\country{1}{France}

	\department{2}{Institut Camille Jordan}
	\organisation{2}{Universit\'e Claude Bernard Lyon 1}
	\rorid{2}{029brtt94}
	\address{2}{43 boulevard du 11 novembre 1918}
	\zip{2}{69622}
	\city{2}{Villeurbanne Cedex}
	\country{2}{France}
	\affemail{2}{nahari@math.univ-lyon1.fr}}
\Emsaffil{2}{
	\department{Institut Camille Jordan}
	\organisation{Universit\'e Claude Bernard Lyon 1}
	\rorid{029brtt94}
	\address{43 boulevard du 11 novembre 1918}
	\zip{69622}
	\city{Villeurbanne Cedex}
	\country{France}
	\affemail{strobl@math.univ-lyon1.fr}}

\classification[57R30]{53C12}


\begin{abstract}
We recall the  notion of  a singular foliation (SF) on a manifold $M$, viewed as an appropriate submodule of $\mathfrak{X}(M)$, and adapt it to the presence of a Riemannian metric $g$, yielding a module version of a singular Riemannian foliation (SRF). Following Garmendia-Zambon on Hausdorff Morita equivalence of SFs, we define the Morita equivalence of SRFs (both in the module sense as well as in the more traditional geometric one of Molino) and show that the  leaf spaces of Morita equivalent SRFs 
are isomrophic as pseudo-metric spaces.

\noindent In a second part, we introduce the category of  $\mathcal{I}$-Poisson manifolds. Its objects and morphisms generalize Poisson manifolds and morphisms in the presence of appropriate ideals $\mathcal{I}$ of the smooth functions on the manifold such that two conditions are satisfied: $(i)$ The category of Poisson manifolds becomes a full subcategory when choosing $\mathcal{I}=0$  and $(ii)$ there is a reduction functor from this new category to the category of Poisson algebras, which generalizes coistropic reduction to the singular setting. 

\noindent Every SF on $M$ gives rise to an  $\mathcal{I}$-Poisson manifold on $T^*M$ and $g$ enhances this to an SRF if and only if the induced Hamiltonian lies in the normalizer of  $\mathcal{I}$. This perspective provides, on the one hand, a simple proof of the  fact that every module SRF is a geometric SRF and, on the other hand, a construction of an algebraic invariant of  singular foliations: Hausdorff Morita equivalent SFs have isomorphic reduced Poisson algebras.
\end{abstract}

\maketitle


\section{Introduction}
\addtocontents{toc}{\protect\setcounter{tocdepth}{1}}

\noindent The first purpose of this article is to introduce and study a notion of singular Riemannian foliations which is adapted to the module definition of a singular foliation. More precisely, following  I. Androulidakis and G. Skandalis \cite{skandalis}, a singular foliation is defined as follows:\footnote{A singular foliation can be equivalently defined as an involutive and locally finitely generated subsheaf of the sheaf of smooth vector fields on $M$ closed under multiplication by $C^\infty(M)$ \cite{Lieinfinity} (see also \cite{garmendiazambon}). This has the advantage  that one can replace $C^\infty(M)$ by an arbitrary sheaf of rings $\mathcal{O}$ on $M$. Definition \ref{def1}, however, is more convenient for the present purposes.}

\begin{Definition}\label{def1} A \emph{singular foliation} (SF) on $M$ is defined as a $C^{\infty}(M)$-submodule $\mathcal{F}$ of the module of compactly supported vector fields on $M$, which is locally finitely generated and closed with respect to the Lie bracket of vector fields.
\end{Definition}


\noindent This definition induces a decomposition of $M$ into injectively immersed submanifolds called leaves \cite{Hermann}, thus yielding singular foliations in the more traditional sense (see, e.g., \cite{lavau}). But the association is not one-to-one: several SFs give rise to the same leaf decomposition. However, in the case where all the leaves have the same dimension, the relation is one-to-one and Definition \ref{def1}  becomes equivalent to the usual notion of a regular foliation. Examples of SFs are induced on the underlying manifold by, e.g., Poisson manifolds, Lie algebroids, and Lie infinity algebroids.

\noindent Now let us add a Riemannian structure $g$ to the above setting. Inspired by \cite{ks16,stroblkotov}, 
but stripping off unnecessary data from the definitions given there, we propose

\begin{Definition}\label{msrf}
A \emph{singular Riemannian foliation} (SRF) on a Riemannian manifold $(M,g)$ is defined as an SF $\mathcal{F}$ on $(M,g)$ such that for every vector field $X\in\mathcal{F}$ we have

\begin{equation} 
    \mathcal{L}_Xg\in \Omega^1(M)\,\odot\,g_{\flat}(\mathcal{F}), \label{LieX}
\end{equation}
where $g_\flat \colon \mathfrak{X}(M)\to \Omega^1(M)$, $X\mapsto g(X,\cdot)$ is the standard musical isomorphism and $\odot$ stands for the symmetric tensor product.
\end{Definition}

\noindent With this definition, every geodesic perpendicular to one leaf turns out to stay perpendicular to all the leaves it meets, thus yielding singular Riemannian foliations in the more traditional sense \cite{molino}. The converse is not always true: A singular Riemannian foliation in the sense of Molino is not always an SRF. For a regular foliation, Definition  \ref{msrf} becomes equivalent to the usual notion of a (regular) Riemannian foliation (\cite{H58}, \cite{R59}). 
Examples of SRFs are given by isometric Lie group actions on Riemannian manifolds and, more generally, orbit decompositions induced by Riemannian groupoids \cite{delhoyo}.

\noindent  Our notion of SRFs behaves well under the \emph{pullback operation} of \cite{skandalis}. This permits us to provide a definition of \emph{Morita equivalence} between SRFs. It implies \emph{Hausdorff Morita equivalence} for the underlying SFs, as defined in \cite{garmendiazambon}. In the fore-cited work it is shown that the leaf spaces of Hausdorff Morita
equivalent SFs are homeomorphic. Here we will establish

\begin{theoremint} \label{theorem1}
Let $(N_1,g_1,\mathcal{F}_1)$ and $(N_2,g_2,\mathcal{F}_2)$ be Morita equivalent SRFs. Then their leaf spaces are isometric as pseudo-metric spaces.
\end{theoremint}

\noindent A second purpose of this article is to introduce the \emph{category of $\mathcal{I}$-Poisson manifolds} \textbf{IPois}. For its objects, the intention is to generalize coisotropic submanifolds (see, e.g.,  \cite{marsdenratiu}) to the singular setting. For simplicity of the 
presentation, in the Introduction we provide the definition of objects for the subcategory of \emph{semi-strict $\mathcal{I}$-Poisson manifolds} \textbf{ssIPois} of \textbf{IPois}, which are constructed simply out of Poisson manifolds:\footnote{For the complete version see Definitions \ref{ipoisson}, \ref{ipoissoncat}, and \ref{semistrict} below. The more general notion permits to cover also examples such as Hamiltonian quasi-Poisson manifolds \cite{anton}, see Example \ref{quasi-poisson}.}

\begin{Definition}\label{def}
A semi-strict \emph{$\mathcal{I}$-Poisson manifold} is a triple $(P,\left\{\cdot,\cdot\right\},\mathcal{I})$ where $\mathcal{I}$ is a subsheaf of smooth functions on a Poisson manifold $(P,\left\{\cdot,\cdot\right\})$ which is closed under multiplication by smooth functions, locally finitely generated, and for every open subset $U\subset P$, $\mathcal{I}(U)\subset C^\infty(U)$ is a Poisson subalgebra, \emph{i.e.}\
\begin{equation*}
\left\{\mathcal{I}(U),\mathcal{I}(U)\right\}\subset \mathcal{I}(U).
\end{equation*}
\end{Definition}

\noindent To describe dynamics, one needs 
a compatible Hamiltonian, i.e.\ a function $H \in N(\mathcal{I})$ where 
\[N(\mathcal{I}):=\left\{f\in C^{\infty}(P)\colon \left\{f\vert_U,\mathcal{I}(U)\right\}\subset \mathcal{I}(U)\,\,\textit{for every open subset }\, U\right\}.\]
We then call $(P,\left\{\cdot,\cdot\right\},\mathcal{I},H)$ a (semi-strict) dynamical $\mathcal{I}$-Poisson manifold and the corresponding category \textbf{(ss)dynIPois}.

\noindent The property that a singular foliation is locally finitely generated is crucial for the existence of the induced leaf decomposition. Similarly, the condition ``locally finitely generated'' in Definition \ref{def} 
is essential for showing  that the flow of any $H \in N(\mathcal{I})$, if complete, preserves the sheaf $\mathcal{I}$ (see Proposition \ref{normcons} for the precise statement). 

\begin{Definition} \label{morph}
A smooth map $\varphi \colon P_1\to P_2$ between   $(P_1,\left\{\cdot,\cdot\right\}_1,\mathcal{I}_1)$ and $(P_2,\left\{\cdot,\cdot\right\}_2,\mathcal{I}_2)$ is a \emph{morphism of (semi-strict) $\mathcal{I}$-Poisson manifolds}, iff the two obvious conditions \\$\varphi^* \left(\mathcal{I}_2(P_2) \right)\subset \mathcal{I}_1(P_1)$ and $\varphi^* N(\mathcal{I}_2) \subset N(\mathcal{I}_1)$ are complemented by

\begin{equation} \label{cond}
 \{\varphi^* f,\varphi^* g\}_1 - \varphi^* \{f,g\}_2 \in  \mathcal{I}_1(P_1) \qquad \forall f,g \in N(\mathcal{I}_2).
\end{equation}
For dynamical $\mathcal{I}$-Poisson manifolds we add the condition $\varphi^* H_2-H_1\in  \mathcal{I}_1$.
\end{Definition}

\noindent These are also the morphisms of the general category, when ``semi-strict'' in the parenthesis is dropped.
With this notion of morphisms, the category \textbf{Pois} of Poisson manifolds is a full subcategory of \textbf{(ss)IPois} for the choice of the zero ideal. In general, however, the morphisms between (semi-strict) $\mathcal{I}$-Poisson manifolds are not necessarily Poisson maps between the underlying Poisson manifolds---an important feature in several applications.

\noindent The condition \eqref{cond} is optimal to ensure that $\varphi^*$ descends to a Poisson morphism on the level of reductions: In fact, every (semi-strict) $\mathcal{I}$-Poisson manifold $(P,\left\{\cdot,\cdot\right\},\mathcal{I})$ induces a Poisson algebra structure on $N(\mathcal{I})/\mathcal{I}(P)$. In the case of coisotropic reductions \cite{marsdenratiu}, this algebra coincides with the algebra of smooth functions on the reduced Poisson manifold. The algebraic formulation here is, however, also applicable in the general context of $\mathcal{I}$-Poisson manifolds, where, e.g., the vanishing set of the ideal $\mathcal{I}$(P) does not need to be a submanifold anymore. The conditions in Definition \ref{morph} ensure that  there is a canonical contravariant functor $F$ from \textbf{(ss)IPois} to \textbf{PoisAlg}, the category of Poisson algebras. 

\noindent The final purpose of this article is to bring the two aforementioned subjects together and, in particular, to use $\mathcal{I}$-Poisson geometry so as to learn more about SFs and SRFs. 

\noindent Starting from an SF $(M,\mathcal{F})$ and viewing every vector field in $\mathcal{F}$ as a smooth function on $T^*M$, we construct a semi-strict $\mathcal{I}$-Poisson manifold $(T^*M,\left\{\cdot,\cdot\right\}_{T^*M},\mathcal{I}_{\mathcal{F}})$. Moreover, every metric $g$ on $M$ defines a compatible Hamiltonian (making the semi-strict $\mathcal{I}$-Poisson manifold dynamical) if and only if the metric satisfies condition \eqref{LieX}. Thus there is a canonical map from SFs and SRFs to the objects of \textbf{ssIPois} and  \textbf{ssdynIPois}, respectively. 
As we will see, this construction is not only conceptually illuminating, it also has technical advantages: we will use it to find elegant proofs of several properties of SFs and SRFs like to show, e.g., that Definition \ref{msrf} automatically induces an SRF in the sense of  \cite{molino}. 

\noindent To complete the above map on objects to a functor, one would need a proper definition of the categories \textbf{SF} and \textbf{SRF} of singular (Riemannian) foliations. Surprisingly, already for SFs,  in the literature there is not yet any satisfactory proposal for what a morphism between general SFs should be. However, the situation changes if one restricts to submersions and Riemannian submersions in the case of SFs and SRFs, respectively, because in these cases the previously mentioned pullback operations are defined. For example, a Riemannian submersion $\pi\colon (N,h)\to (M,g)$ between two SRFs $(N,h,\mathcal{F}_N)$ and $(M,g,\mathcal{F}_M)$
which satisfies $\pi^{-1}\mathcal{F}_M=\mathcal{F}_N$ should definitely be considered as a morphism. Let us call \textbf{SF}$_0$ and \textbf{SRF}$_0$
the two (sub)categories with such restricted morphisms. In this paper we show in particular

\begin{theoremint} \label{thmfunctor}
There are canonical functors $\Psi\colon \mathbf{SF}_0 \to \mathbf{IPois}$ and $\Phi \colon \mathbf{SRF}_0 \to \mathbf{dynIPois}$.
\end{theoremint}

\noindent As a side result, we will find that for $\mathcal{F}_M=0$, $\Phi(\pi)$ becomes an ordinary Poisson map if and only if the horizontal distribution $\left(\ker \exd \pi\right)^\perp$ is integrable---correcting \cite{baez}, where this map has been considered as well, but claimed to always be Poisson. 

\noindent Composing the functor $\Psi$, evaluated on an SF $(M,\mathcal{F})$, with the functor $F \colon \textbf{IPois} \to \textbf{PoisAlg}$, we obtain the (reduced) Poisson algebra $A(\mathcal{F}) : = N(\mathcal{I}_\mathcal{F})/\mathcal{I}_\mathcal{F}(T^*M)$. This algebra provides an invariant of Hausdorff Morita equivalence, since we will prove
 
\begin{theoremint}\label{Moritainv}
Let $(M_1,\mathcal{F}_1)$ and $(M_1,\mathcal{F}_1)$ be Hausdorff Morita equivalent singular foliations. Then the reduced Poisson algebras  $A(\mathcal{F}_1)$ and $A(\mathcal{F}_2)$ are isomorphic.
\end{theoremint}

\noindent The structure of this paper is as follows:

\noindent \emph{Section \ref{sfrev}} contains a short review of the definitions and main properties of SFs related to the goal of this paper, in particular the notion of Hausdorff Morita equivalence of SFs. 

\noindent In \emph{Section \ref{srfsection} } we introduce SRFs and study some of their properties. 
We show (in Theorem \ref{thm:almostKilling} below)
that every finitely generated SRF  
admits  an almost Lie algebroid structure with connection to turn the SRF into an almost Killing  Lie algebroid \cite{stroblkotov}. We  define Morita equivalence of SRFs, show that it defines an 
equivalence relation, and prove Theorem \ref{theorem1}.

\noindent \emph{Section \ref{i-poisson}} introduces the category \textbf{IPois}, the reduction functor $F$ to \textbf{PoisAlg}, and provides several examples and properties of $\mathcal{I}$-Poisson manifolds.

\noindent In \emph{Section \ref{objects}} we show how SFs and SRFs give rise to particular  $\mathcal{I}$-Poisson and dynamical $\mathcal{I}$-Poisson manifolds, respectively.

\noindent In
\emph{Section \ref{functor}}, finally, we prove Theorems \ref{thmfunctor} and \ref{Moritainv}.

\noindent The definition of almost Killing Lie algebroids as well as part of the proof of Theorem \ref{thm:almostKilling} (in the form of Proposition \ref{prop:almostKilling}) are deferred to \emph{Appendix
\ref{app}}.

\section{Background on singular foliations and their Morita equivalence}\label{sfrev}

\noindent In what follows, $M$ is assumed to be a smooth manifold and $\mathfrak{X}_c(M)$ denotes the $C^{\infty}(M)$-module of compactly supported vector fields on $M$. For more details and examples of singular foliations see \cite{skandalis} or \cite{Lieinfinity}. 

\begin{Definition}\label{lfg}
A $C^{\infty}(M)$-submodule $\mathcal{F}\subset\mathfrak{X}_c(M)$ is called \emph{locally finitely generated} if for every point $q\in M$ there exist an open neighborhood $U\subset M$ around $q$ such that the submodule $\iota_U^{-1}\mathcal{F}\subset \mathfrak{X}_c(U)$ defined as
\[\iota_U^{-1}\mathcal{F}:=\left\{X\in \mathcal{F}\,\colon\,\mathrm{supp}(X)\subset U\right\}\]
is finitely generated; \emph{i.e.}\ there exist finitely many vector fields $X_1,\ldots,X_N\subset \mathfrak{X}(U)$ for some positive integer $N$, such that 
\[\iota_U^{-1}\mathcal{F}=\la X_1,\ldots,X_N\ra_{C^\infty_c(U)}.\]
\end{Definition}

\begin{Remark}
Note that the generators of $\iota_U^{-1}(\mathcal{F})$ in Definition \ref{lfg} are not required to be compactly supported. This allows for more flexibility in constructing singular foliations on $M$.
\end{Remark}

\begin{Example} \label{cexpl}
Let $M=\R$. Then the $C^\infty(\R)$-module $\mathfrak{X}_c(\R)$ is globally generated by the single vector field 
$\tfrac{\exd}{\exd x}$. On the other hand, the $C^{\infty}(\R)$-submodule of compactly supported vector fields which vanish on $\R_{-}$ is not locally finitely generated around $0$. 
\end{Example}

\begin{Definition} A \emph{singular foliation} on $M$---SF for short---is defined as $C^{\infty}(M)$-submodule $\mathcal{F}$ of $\mathfrak{X}_c(M)$ which is locally finitely generated and closed with respect to the Lie bracket of vector fields. The pair $(M,\mathcal{F})$ is then called a \emph{foliated manifold}. \label{SF}
\end{Definition}

\begin{Remark}
 One can equivalently define SFs as an involutive and locally finitely generated subsheaf of the sheaf of vector fields $\mathfrak{X}$. This is equivalent to Definition \ref{SF} in the smooth setting, but it has advantages if we wish to work with the sheaves of algebraic, real analytic or holomorphic functions (See \cite{Lieinfinity} or \cite{garmendiazambon}). In particular, since the mentioned sheaves of rings are Noetherian, the condition of being locally finitely generated is automatically satisfied and therefore can be dropped.
\end{Remark}

\noindent A classical theorem of R. Hermann \cite{Hermann} implies that an SF defined as above partitions $M$ into smooth, connected, and injectively immersed submanifolds (of possibly different dimensions) called \emph{leaves}.

\noindent Let $L_q$ be the leaf passing through the point $q\in M$ in a foliated manifold $(M,\mathcal{F})$. Then, by definition of the leaves, $T_qL_q$ can be identified with $\left\{X\vert_q\,\colon\,X\in \mathcal{F}\right\}\subset T_qM$, which motivates

\begin{Definition}
For every point $q\in M$ in a foliated manifold $(M,\mathcal{F})$, the \emph{tangent of $\mathcal{F}$} at $q$ is defined as

\begin{equation*}
    F_q:=\left\{X\vert_q\,\colon\, X\in \mathcal{F}\right\}\subset T_qM.
\end{equation*}
\end{Definition}

\noindent If $q \mapsto \dim(F_q)$ is constant on $M$, we obtain \emph{regular foliations} as particular singular ones. In this case, by the Frobenius theorem, there is a one-to-one correspondence between the leaf decomposition of the foliation and the module of vector fields generating it. This is no more the case if the SF is non-regular; there \emph{always} exist different modules generating the same leaf decomposition then (for an example, see Example \ref{k} below). Note also that in the non-regular case all the vector fields tangent to the leaves of a given SF do not necessarily define an SF anymore: the module of Example \ref{cexpl}, despite not being an SF since not finitely generated, induces a leaf decomposition, which can be obtained also by an honest SF with the single generator $\chi \,  \tfrac{\exd}{\exd x}$. Here $\chi \in C^\infty(\R)$ can be chosen, e.g., as the  function

\begin{equation}
\label{chi}
\chi(x) =  \begin{cases}
\exp(\tfrac{-1}{x^2}) &\quad x>0\\
0 &\quad x\leq 0
\end{cases}\,\cdot
\end{equation}  

\begin{Remark}
The function $M\to \N$ given by $q\mapsto \dim(F_q)$ is lower semi-continuous. As a result, the subset $U\subset M$ of the continuity set of $\dim(F_q)$ is open and dense, and $\mathcal{F}|_U$ induces a regular foliation over each connected component of $U$ \cite{skandalis}.
\end{Remark}

\noindent The following example shows the importance of being locally finitely generated.

\begin{Example}
On $M=\R^2$, consider the module $\mathcal{G}$
  generated by the vector fields $\partial_{x}$ and $\mathcal{X}(x)\partial_{y}$, where $\chi$ is the function defined in \eqref{chi}, together with all their multiple commutator Lie brackets. Then, by construction, $\mathcal{G}$ is closed under the Lie bracket. However, it is not locally finitely generated as a $C^\infty(M)$-module since with each derivative on $\mathcal{X}$ we obtain a new, independent coefficient in front of $\partial_{y}$. As a consequence, we loose the well-behavedness of a leaf-decomposition: although every two points in $\R^2$ can be connected by a sequence of flows of vector fields in $\mathcal{G}$ (so that, in this sense, there would be only one leaf that is $\R^2$ itself), the tangent of $\mathcal{G}$ at every point in the left half-plane is only one-dimensional.
\end{Example}

\noindent As mentioned above, Definition \ref{SF} contains more information than a well-behaved decomposition of $M$ into leaves.

\begin{Example} \label{k}
Let $M = \R$ and let $\mathcal{F}$ be an SF generated by vector fields vanishing at the origin of at least order $k \in \N$. The leaf decomposition induced by $\mathcal{F}$ is $\R = \R_- \cup {0} \cup \R_+$ for every choice of $k$. Thus, $\mathcal{F}$ is not  completely determined  by its leaf decomposition.
\end{Example}

\noindent To capture some of this additional information contained in the definition of an SF, we extract some more data from the module $\mathcal{F}$ by the following definition of \cite{skandalis}.

\begin{Definition}
Let $(M,\mathcal{F})$ be a foliated manifold. For every point $q\in M$, the \emph{fiber of $\mathcal{F}$} at $q$ is defined as:

\begin{equation*}
    \mathcal{F}_q:=\mathcal{F}/I_q\!\cdot\!\mathcal{F}
\end{equation*}
where $I_q:=\left\{f\in C^{\infty}(M)\,\colon\,f(q)=0\right\}$ is the vanishing ideal of $q$ in $C^{\infty}(M)$.
\end{Definition}

\begin{Remark}
The function $M\to \N$ sending $q\to \dim(\mathcal{F}_q)$ is upper semi-continuous, and $\dim(\mathcal{F}_q)$ gives the minimal number of vector fields locally generating $\mathcal{F}$ around $q\in M$ \cite{skandalis}.
\end{Remark}

\noindent Note that for every point $q\in M$ the evaluation map $ev_q\colon \mathcal{F}_q\to F_q$, $[X]\to X\vert_q$ is a homomorphism of vector spaces and induces the following short exact sequence:
\[0\to \ker(ev_q)\to \mathcal{F}_q\to F_q\to 0\]

\noindent here $[X]$ denotes the equivalence class of the vector field $X\in \mathcal{F}$.

\noindent It is not difficult to see that the Lie bracket on $\mathcal{F}$ induces a Lie bracket on the finite-dimensional vector space $\ker(ev_q)\subset \mathcal{F}_q$.

\begin{Definition}
The vector space $\mathfrak{g}_q^{\mathcal{F}}:=\ker(ev_q)$ together with the bracket inherited by $\mathcal{F}_q$ defines the \emph{isotropy Lie algebra of $\mathcal{F}$} at $q$.
\end{Definition}

\noindent In the case of regular foliations, the map $ev_q\colon \mathcal{F}_q\to F_q$ is a vector space isomorphism and $\mathfrak{g}_q^{\mathcal{F}}=0$. So one can say that the isotropy Lie algebra $\mathfrak{g}_q^{\mathcal{F}}$ characterizes in part the singularity of $\mathcal{F}$ at $q\in M$.

\noindent In Example \ref{k} all fibers and isotropy Lie algebras at the origin are isomorphic. This changes, if we increase the dimension of $M$.

\begin{Example}
Let $M=\R^n$, $n\geq 2$, and let $\mathcal{F}$ be the SF generated by vector fields vanishing at the origin at least of order $k\in\N$. There are always only  two leaves $M\backslash \left\{0\right\}$ and $\left\{0\right\}$, but the fiber at the origin has different dimensions for different choices of $k$,  $\dim \mathcal{F}_0=\binom{k+n-1}{n-1}$. 
\end{Example}

\begin{Definition} \label{defpullback} Let $(M,\mathcal{F})$ be an SF and $\pi \colon N \to M$ a submersion, then the $C^\infty(N)$-module generated by vector fields on $N$ projectable to $\mathcal{F}$ defines the \emph{pullback foliation} $\left(N,\pi^{-1}\mathcal{F}\right)$. 
\end{Definition}

\noindent Here a vector field $V$ on $N$ is called projectable to $\mathcal{F}$ if there exists a vector field $X\in \mathcal{F}$ such that for every point $q\in N$ we have

\begin{equation*}
    \exd _q\pi(V\vert_q)=X\vert_{\pi(q)}.
\end{equation*}
As shown in  \cite{skandalis}, Propositions $1.10$ and $1.11$, the pullback foliation is indeed finitely generated and involutive, i.e.\ it is an SF. This notion behaves well under composition of submersions: For submersions $\pi_P\colon P\to M$ and $\pi_M \colon M\to N$, one has \begin{equation*}
    (\pi_M\circ \pi_P)^{-1}\mathcal{F}=\pi_P^{-1}(\pi_M^{-1}\mathcal{F}).
\end{equation*}
As an example, if $U$ is an open subset of a foliated manifold $(M,\mathcal{F})$, then for the inclusion map $\iota_U\colon U\hookrightarrow M$, the SF $\iota_U^{-1}\mathcal{F}$ is compatible with Definition \ref{lfg}.

\begin{Definition}[\cite{garmendiazambon}]
Two foliated manifolds $(M_1,\mathcal{F}_1)$ and $(M_2,\mathcal{F}_2)$ are \emph{Hausdorff Morita equivalent} if there exists a smooth manifold $N$ and surjective submersions with connected fibers $\pi_i \colon N\to M_i$, $i=1,2$ such that

\begin{equation*}
    \pi_1^{-1}\mathcal{F}_1=\pi_2^{-1}\mathcal{F}_2.
\end{equation*}

\noindent In this case we write $(M_1,\mathcal{F}_1)\sim_{ME}(M_2,\mathcal{F}_2)$.
\end{Definition}

\noindent It is shown in  \cite{garmendiazambon} that the SFs underlying Morita equivalent Lie algebroids \cite{gin} or  Morita equivalent Poisson manifolds \cite{Xu} are Hausdorff Morita equivalent. Also the Morita equivalence of regular foliations \cite{molino} is a special case. Hausdorff Morita equivalence defines an equivalence relation on foliated manifolds---something that holds true for Poisson manifolds only upon restriction to those integrating to a symplectic groupoid. The main fact about Hausdorff Morita equivalent foliated manifolds is that they have Morita equivalent Holonomy groupoids (as open topological groupoids) defined in \cite{skandalis}.

\begin{theorem}[\cite{garmendiazambon}]\label{hme} Let $(M_1,\mathcal{F}_1)$ and $(M_2,\mathcal{F}_2)$  be foliated manifolds which are Hausdorff Morita equivalent by means of $(N,\pi_1,\pi_2)$. Then
\noindent(i) The map sending the leaf passing through $q\in M_1$ to the leaf of $\mathcal{F}_2$ containing $\pi_2(\pi_1^{-1}(q))$ is a homeomorphism between the leaf spaces. It preserves the codimension of leaves and the property of being an embedded leaf.
\noindent(ii) Let $q_1\in N_1$ and $q_2\in N_2$ be points in corresponding leaves. Choose transversal slices $S_{q_1}$ at $q_1$ and $S_{q_2}$ at $q_2$. Then the foliated manifolds $(S_{q_1},\iota_{S_{q_1}}^{-1}\mathcal{F}_1)$ and $(S_{q_2},\iota_{S_{q_2}}^{-1}\mathcal{F}_2)$ are diffeomorphic and the isotropy Lie algebras $\mathfrak{g}_{q_1}^{\mathcal{F}_1}$ and $\mathfrak{g}_{q_2}^{\mathcal{F}_2}$ are isomorphic.
\end{theorem}

\begin{Example}
 For smooth, connected manifolds $M$ and $N$, $(M,\mathfrak{X}_c(M))$ and $(N,\mathfrak{X}_c(N))$ are always Hausdorff Morita equivalent. 
 On the other hand, $(M,0)$ and $(N,0)$ are Hausdorff Morita equivalent only if $M$ and $N$ are diffeomorphic.
\end{Example}

\section{Singular Riemannian foliations and their Morita equivalence}\label{srfsection}

\noindent In what follows, $(M,g)$ denotes a Riemannian manifold. We first recall the traditional notion of a \emph{singular Riemannian foliations} (SRF) motivated by \cite{molino}, to which we will add the suffix ``geometric'' so as to distinguish it from a second one that we will introduce directly below.

\begin{Definition}\label{gsrf}
Let $\mathcal{F}$ be an SF on $(M,g)$. We call the triple $(M,g,\mathcal{F})$ a \emph{geometric SRF}, if every geodesic orthogonal to a leaf at one point is orthogonal to all the leaves it meets.
\end{Definition}

\noindent In this text, we focus mainly on the following definition of SRFs, streamlining the one given in \cite{stroblkotov}\footnote{For the relation of module SRFs with the notion defined in \cite{stroblkotov} see Appendix \ref{app} as well as Theorem \ref{thm:almostKilling} below.}:

\begin{Definition}\label{aSRF}
Let $\mathcal{F}$ be an SF on $(M,g)$. We call the triple $(M,g,\mathcal{F})$  a \emph{module SRF}, if for every vector field $X\in\mathcal{F}$ we have

\begin{equation} \label{asrfcond}
    \mathcal{L}_Xg\in \Omega^1(M)\,\odot\,g_{\flat}(\mathcal{F}).
\end{equation}
\end{Definition}
\noindent Here $\odot$ stands for the symmetric tensor product and $g_{\flat}$ is the map on sections induced by the musical isomorphism $g_{\flat}\colon TM\to T^*M$, $(q,v)\mapsto g_q(v,\cdot)$. Let $(g_{\flat})^{-1} \colon \Omega^1(M) \to \mathcal{X}(M)$ denote the corresponding inverse map and $g^{-1} \in \Gamma(S^2TM)$ the 2-tensor inducing it. Then, by means of $\mathcal{L}_X(g_{\flat})^{-1}=-(g_{\flat})^{-1} \circ (\mathcal{L}_Xg_{\flat}) \circ (g_{\flat})^{-1}$, we can express the defining property of a module SRF also in the following form

\begin{Lemma}\label{lem1}
The triple $(M,g,\mathcal{F})$ is a module SRF if and only if
\[\mathcal{L}_Xg^{-1}\in \mathfrak{X}(M)\odot \mathcal{F}\]
 for every vector field $X\in \mathcal{F}$.
\end{Lemma}

\noindent As a consequence of the following lemma and proposition, it is enough to check Equation \eqref{asrfcond} locally for a family of generators.

\begin{Lemma} \label{srfgen}
Let $(M,\mathcal{F})$ be a foliated manifold such that $\mathcal{F}=\left\la X_1,\ldots,X_N\right\ra_{C^{\infty}_c(M)}$ for some positive integer $N$. Then the triple $(M,g,\mathcal{F})$ is a module SRF if and only if there exist $\omega_a^b\in\Omega^1(M)$ for $a,b=1,\ldots,N$ such that
\begin{equation*}
   \mathcal{L}_{X_a}g=\sum_{b=1}^N\omega_a^b\,\odot\,g_{\flat}(X_b).
\end{equation*}
\end{Lemma}

\noindent \begin{proof}
First assume that $(M,g,\mathcal{F})$ is a module SRF. Choose a partition of unity $\left\{\rho_i\right\}_{i=1}^\infty$ subbordinate to a locally finite cover $\left\{U_i\right\}_{i=1}^\infty$ of $M$. For every $a=1,\ldots,N$ we have
\begin{align*}
    \mathcal{L}_{X_a}g=\sum_{i=1}^\infty\rho_i\mathcal{L}_{X_a}g&=\sum_{i=1}^\infty\left(\mathcal{L}_{\rho_iX_a}g-(\exd\rho_i)\,\odot\,g_\flat(X_a)\right)\\
    &=\sum_{i=1}^\infty\left(\sum_{b=1}^N\eta_{i,a}^b\,\odot\,g_\flat(X_b)-(\exd\rho_i)\,\odot\,g_\flat(X_a)\right)\\
    &=\sum_{b=1}^N\omega_a^b\,\odot\,g_{\flat}(X_b),
\end{align*}
for some $1$-forms $\eta_{i,a}^b$ on $M$ and $\omega_a^b:=\sum_{i=1}^\infty\eta_{i,a}^b-\delta_a^b\exd\rho_i$.
For the converse,  let $X$ be a vector field in $\mathcal{F}$. By assumption, there exist $f^1,\ldots,f^N\in C^{\infty}_c(M)$ such that $X=\sum_{a=1}^Nf^aX_a$. It follows that

\begin{align*}
\mathcal{L}_Xg=\sum_{a=1}^N\mathcal{L}_{f^aX_a}g&=\sum_{a=1}^Nf^a\mathcal{L}_{X_a}g+ (\exd f^a)\,\odot\,g_{\flat}(X_a)\\
&=\sum_{a=1}^N\omega_a^b\,\odot\,g_{\flat}(f^aX_b)+(\exd f^a)\,\odot\,g_{\flat}(X_a)\in \Omega^1(M)\,\odot\,g_{\flat}(\mathcal{F}),
\end{align*}
\end{proof}

\noindent An important property of the definition of a geometric SRF is that the defining condition is local. This is less trivial in the case of module SRFs.

\begin{Proposition} \label{locality}
The triple $(M,g,\mathcal{F})$ is a module SRF if and only if for every point $q\in M$ there exist an open neighborhood $U\subset M$ around $q$ such that $(U,g_U,\iota_U^{-1}\mathcal{F})$ is a module SRF, where $g_U$ is the restriction to $U$ of $g$.
\end{Proposition}

\noindent \begin{proof}
If $(M,g,\mathcal{F})$ is a module SRF, then restricting both sides of Equation \eqref{asrfcond} to any open subset $U\in M$ implies that $(U,g_U,\iota_U^{-1}\mathcal{F})$ is a module SRF. It remains to prove the converse. Choose a partition of unity $\left\{\rho_i\right\}_{i=1}^\infty$ subbordinate to a locally finite cover $\left\{U_i\right\}_{i=1}^\infty$ of $M$, with open subsets $U_i$ small enough such that $\iota_{U_i}^{-1}\mathcal{F}=\left\la X_{i,1},\ldots,X_{i,N_i}\right\ra_{C^{\infty}_c(U_a)}$ for some positive integer $N_i$ and vector fields $X_{i,1},\ldots,X_{i,N_i}\in\mathfrak{X}(U_i)$. Then for every vector field $X\in\mathcal{F}$,
\[X=\sum_{i=1}^{\infty}\rho_iX.\]

\noindent Moreover, for every positive integer $i$, there exist functions $f^{i,1},...,f^{i,N_i}\in C^{\infty}_c(U_i)$ such that
\[\rho_iX
=\sum_{a=1}^{N_i}f^{i,a}X_{i,a},\]
and consequently 
\[X=\sum_{i=1}^{\infty}\sum_{a=1}^{N_i}f^{i,a}X_{i,a}.\]
This together with Lemma \ref{srfgen} now permit us to prove that $(M,g,\mathcal{F})$ is a module SRF. We have

\begin{align*}
\mathcal{L}_Xg=\sum_{i=1}^{\infty}\sum_{a=1}^{N_i}f^{i,a}\mathcal{L}_{X_{i,a}}g_{U_i}+(\exd f^{i,a})\,\odot\,(g)_{\flat}(X_{i,a}),
\end{align*}
which proves $\mathcal{L}_Xg\in \Omega^1(M)\,\odot\,g_{\flat}(\mathcal{F})$ since $X$ is compactly supported and only finitely many $f^{i,a}$ are nonzero on $\mathrm{supp}(X)$.
\end{proof}

\noindent Every finitely generated SF is image of the anchor map of an almost Lie algebroid  \cite{Lieinfinity} (see Appendix \ref{app}). For module SRFs, one has furthermore

\begin{theorem}\label{thm:almostKilling} Let $(M, g)$ be a Riemannian manifold. The following statements hold true:
\begin{enumerate}[(i)]
    \item For every module SRF $(M,g,\mathcal{F})$ with  $\mathcal{F}$ finitely generated, there exists an almost Lie algebroid $(A,\rho,[\cdot,\cdot]_A)$ over $M$ equipped with a connection $\nabla\colon\Gamma(A)\to\Gamma(T^*M\otimes A)$ such that $\mathcal{F}:=\rho(\Gamma_c(A))$ and 
    
    \begin{equation}
        {}^A\nabla g = 0, \label{Anablag}
    \end{equation}
    where ${}^A\nabla$ is the $A$-connection induced by $\nabla$, see  Equation \eqref{tau} in Appendix  \ref{app}.
    \item Let $(A,\rho,[\cdot,\cdot]_A)$ be an almost Lie algebroid over a Riemannian manifold $(M,g)$, such that the triple $(M,g,\mathcal{F}:=\rho(\Gamma_c(A)))$ is a module SRF. Then there exists a connection $\nabla$ on $A$ such that \eqref{Anablag} holds true.
\end{enumerate}
\end{theorem}

\noindent\begin{proof}
The proof of the first part of the Theorem can be performed by a straightforward adaptation of the proof of  Proposition \ref{prop:almostKilling} in the Appendix. In particular, the almost Lie algebroid $A$ then can be chosen to be trivial, $A=M \times \R^r$, where $r$ is the number of generators of $\mathcal{F}$. 

\noindent We prove the second  part of the Theorem, where now one is given a particular, not necessarily trivial almost Lie algebroid $A$ inducing $\mathcal{F}$, as follows: There exists a vector bundle $V\to M$ such that $(\Tilde{A}:=A\oplus V)\to M$ is a trivial vector bundle of rank $N$. Consequently there exist sections $e_1,\ldots,e_N\in \Gamma(A)$ and $v_1,\ldots,v_N\in \Gamma(V)$ such that $e_1+v_1,\ldots,e_N+v_N$ is a global frame for $\Tilde{A}$. Now we define the almost Lie algebroid $(\Tilde{A},\Tilde{\rho},[\cdot,\cdot]_{\Tilde{A}})$, where the bracket and the anchor map are the trivial prolongation of $[\cdot,\cdot]_A$ and $\rho$ to $\Tilde{A}$ (since in an almost Lie algebroid one does not need to satisfy the Jacobi identity for the bracket, this extension does not pose any problems here). By assumption $\Tilde{\rho}(\Gamma_c(\Tilde{A}))=\rho(\Gamma_c(A))$ defines a module SRF on $(M,g)$. According to Lemma \ref{srfgen}, this is equivalent to the existence of $1$-forms $\omega_a^b\in \Omega^1(M)$ such that

\begin{equation}
    \mathcal{L}_{X_a}g=\sum_{b=1}^N\omega_a^b\odot \iota_{X_b}g\,\,\,\,\,\forall\,a=1,\ldots,N.
\end{equation}
Here $X_a:=\Tilde{\rho}(e_a+v_a)=\rho(e_a)$. Now define a connection $\Tilde{\nabla}$ on $\Tilde{A}$ by

\begin{equation*}
    \Tilde{\nabla} (e_a+v_a)=\sum_{b=1}^N\omega_a^b\otimes (e_b+v_b),
\end{equation*}
which induces a connection on $A$ as follows: Let $s\in \Gamma(A)\subset \Gamma(\Tilde{A})$, then 

\begin{equation*}
    \nabla_X s: =\mathrm{Pr}_{A}\circ (\Tilde{\nabla}_X s)\,\, \forall\, X\in\mathfrak{X}(M),
\end{equation*}
where $\mathrm{Pr}_{A} \colon \tilde{A} \to A$ is the projection to the first component.
In particular, for every $e_a$, there exist unique functions $f_a^b\in C^{\infty}(M)$ for $b=1,\ldots,N$ such that $e_a=\sum_{b=1}^Nf_a^b(e_b+v_b)$ and we have

\begin{equation*}
    \nabla_X e_a=Pr_{A}\circ (\Tilde{\nabla}_X (\sum_{b=1}^Nf_a^b(e_b+v_b)))=\sum_{b=1}^NX(f_a^b)e_b +\sum_{b,c=1}^N(f_a^b\iota_X\omega_b^c)e_c.
\end{equation*}
Now for every vector field $X\in \mathfrak{X}(M)$, we have
\begin{align*}
    2g\left(\rho\left(\nabla_Xe_a\right),X\right)&=2g\left(\rho\left(\sum_{b=1}^NX(f_a^b)e_b +\sum_{b,c=1}^N(f_a^b\iota_X\omega_b^c)e_c\right),X\right)\\
                           &=2\sum_{b=1}^Nf_a^b((\sum_{c=1}^N \left(\iota_X\omega_b^c\right)g\left(X_c,X\right))+2\sum_{b=1}^NX(f_a^b)g(X_b,X)\\
                           &=\sum_{b=1}^Nf_a^b\left(\sum_{b=1}^N\omega_b^c\odot \iota_{X_c}g\right)(X,X)+\sum_{b=1}^N(\exd f_a^b\odot \iota_{X_b}g)(X,X)\\
                           &=\sum_{b=1}^N\left(f_a^b\mathcal{L}_{X_b}g+\exd f_a^b\odot\iota_{X_b}g \right)(X,X)\\
                           &=\left(\mathcal{L}_{X_a}g\right)(X,X),
\end{align*}
and, by Lemma \ref{lietau} in the Appendix below, the statement then follows.
\end{proof}

\noindent So locally one can define SFs also as an equivalence class of almost Lie algebroids and module SRFs as an equivalence class of almost Lie algebroids over a Riemannian base with an appropriately compatible connection. (For some related cohomology see also \cite{Sasha-Thomas}).

\noindent Using the language of almost Lie algebroids, the following proposition is Theorem 7 in \cite{stroblkotov}. It will be proven in an alternative, more direct way in the present paper, using the techniques of $\mathcal{I}$-Poisson geometry: 

\begin{Proposition}\label{alggeo}
Every module SRF is a geometric SRF.
\end{Proposition}

\noindent Note that the converse is not true,  at least not for every choice of the module $\mathcal{F}$.

\begin{Example}
    Consider $\mathcal{F}=\la(x^2+y^2)(x\partial_y-y\partial_x)\ra_{C^{\infty}_c(\R^2)}$ on $M=\R^2$ equipped with the standard metric $\exd s^2$. The leaves are circles centered at the origin, which is a geometric SRF, but it does not satisfy Equation \eqref{asrfcond}. More precisely, for $V:=(x^2+y^2)(x\partial_y-y\partial_x)$, a simple calculation implies that
    
    \begin{equation*}
        \cL_V\exd s^2=4\left[\frac{x\exd x+y\exd y}{x^2+y^2}\right]\odot (\exd s^2)_\flat (V),
    \end{equation*}
    on $\R^2\setminus {(0,0)}$. Evidently, the $1$-form $\frac{x\exd x+y\exd y}{x^2+y^2}$ fails to have a smooth extension to the origin.
\end{Example}

\noindent \begin{Remark}
 One can pose the following question as well: Assume that a leaf decomposition of a Riemannian manifold is given, such that the compatibility condition of Definition \ref{gsrf} is satisfied. Is there an SF generating a module SRF with the given leaf decomposition? A counter-example for the polynomial or analytic setting is the singular octonionic Hopf foliation \cite{sohf}: albeit there do exist such (real analytic or polynomial) SFs generating the leaf decomposition, the condition \eqref{asrfcond} is not satisfied for \emph{any} of them. For the smooth setting, this is still an open problem.
\end{Remark}

\noindent For SFs there is a pullback under submersions, see Definition \ref{defpullback} and the text following it. To adapt this to the context of SRFs, we consider the following: 

\begin{Definition}
Let $\pi\colon(N,h)\to(M,g)$ be a smooth submersion between Riemannian manifolds. It is called a \emph{Riemannian submersion} if, for every $q\in N$, the restriction $\exd_q\pi\colon\mathcal{H}_q\to T_{\pi(q)}M$ of $\exd_q\pi$ to $\mathcal{H}_q=(\ker \exd_q\pi)^{\perp_g}\subset T_qN$   is an isometry. The smooth distribution $\mathcal{H}=(\mathcal{H}_q)_{q\in N}$ of rank $\dim(M)$ is called the \emph{horizontal distribution} of $\pi$.
\end{Definition}

\begin{Lemma}\label{decomposition}
Let $\pi \colon (N,h)\to (M,g)$ be a Riemannian submersion and $(M,\mathcal{F})$ an SF. Then the pullback SF can be generated as follows

\begin{equation}
    \pi^{-1}\mathcal{F} = \la\mathcal{F}^{\mathcal{H}}+\Gamma\left(\ker \exd\pi\right)\ra_{C^{\infty}_c(N)},
\end{equation}
where $\mathcal{F}^{\mathcal{H}}$ is the horizontal lift of $\mathcal{F}$.
\end{Lemma}

\noindent \begin{proof}
By definition \ref{defpullback}  the inclusion $\la\mathcal{F}^{\mathcal{H}}+\Gamma\left(\ker \exd\pi\right)\ra_{C^{\infty}_c(N)}\subset \pi^{-1}\mathcal{F}$ is evident. Now let $W$ be a projectable vector field on $N$ projecting to $\mathcal{F}$, i.e.\ there exists a vector field $X\in \mathcal{F}$ such that $\exd_q\pi(W\vert_q)=X\vert_{\pi(q)}$. On the other hand, if we decompose $W$ into its horizontal part $W_H$ and its vertical part $W_V$, we have $\exd_q\pi(W_H\vert_q)=X\vert_{\pi(q)}$, which gives $X_H=V^{\mathcal{H}}$. This means that generators of $\pi^{-1}\mathcal{F}$ belongs to $\mathcal{F}^{\mathcal{H}}+\Gamma\left(\ker \exd\pi\right)$, consequently $\pi^{-1}\mathcal{F} = \la\mathcal{F}^{\mathcal{H}}+\Gamma\left(\ker \exd\pi\right)\ra_{C^{\infty}_c(N)}$.
\end{proof}

\begin{Proposition}\label{pullback}
Let $\pi\colon (N,h)\to (M,g)$ be a Riemannian submersion and let $(M,g,\mathcal{F})$ be a module SRF. Then $(N,h,\pi^{-1}\mathcal{F})$ is a module SRF as well. The same statement holds true for geometric SRFs.
\end{Proposition}

\noindent Proposition \ref{pullback} will be proven in Section \ref{functor} below. As a consequence, and by the fact that (regular) Riemannian foliations are locally modeled on Riemannian submersions \cite{molino}, we obtain

\begin{Proposition}
Let $(M,\mathcal{F})$ be a regular foliation on a Riemannian manifold $(M,g)$. Then $(M,g,\mathcal{F})$ is a geometric SRF if and only if it is a module SRF. 
\end{Proposition}

\begin{Example}
Let $G$ be a Lie group acting by isometries on $(M,g)$. Then after Lemma \ref{srfgen} the $C^{\infty}(M)$-submodule $\mathcal{F}\subset\mathfrak{X}_c(M)$ generated by fundamental vector fields is a module SRF on $(M,g)$, since every fundamental vector field $X$ is a Killing vector field: $\mathcal{L}_Xg=0$.
\end{Example}

\begin{Example}
The proof of Theorem 1 in \cite{stroblkotov} shows that the geometric SRF induced on the manifold of objects of a Riemannian Groupoid---as defined in \cite{delhoyo}---is a module SRF.
\end{Example}

\begin{Definition} \label{riemmorita}
Two module SRFs $(M_1,g_1,\mathcal{F}_1)$ and $(M_2,g_2,\mathcal{F}_2)$ are \emph{Morita equivalent} if there exists a Riemannian manifold $(N,h)$ together with two surjective Riemannian submersions with connected fibers $\pi_i\colon (N,h)\to (M_i,g_i)$ for $i=1,2$ such that

\begin{equation*}
    \pi_1^{-1}\mathcal{F}_1=\pi_2^{-1}\mathcal{F}_2
\end{equation*}
\noindent and we write $(N_1,g_1,\mathcal{F}_1)\sim_{ME}(N_2,g_2,\mathcal{F}_2)$.
\end{Definition}

\begin{Remark} \label{underlying}
This notion of Morita equivalence can be defined for geometric SRFs as well as for module ones. Consequently, if two module SRFs are Morita equivalent then they are also Morita equivalent as geometric SRFs. Moreover, if we forget about Riemannian metrics, we obtain Hausdorff Morita equivalent foliated manifolds.
\end{Remark}

\noindent While for Hausdorff Morita equivalence of SFs transitivity of the equivalence relation is relatively easy to show, this is more involved in case of the additional Riemannian structure due to the presence of the metric.

\begin{Proposition}
The Morita equivalence of module SRFs defines an equivalence relation. 
\end{Proposition}

\noindent \begin{proof}
Reflexivity is evident from the definition and for the self-equivalence the identity map defines a Morita equivalence between a module SRF and itself. Now we prove the transitivity as follows: Assume that $(M_1,g_1,\mathcal{F}_1)\sim_{ME}(M_2,g_2,\mathcal{F}_2)$ given by $\pi_i\colon (U,g_U)\to (M_i,g_i)$ for $i=1,2$ and $(M_2,g_2,\mathcal{F}_2)\sim_{ME}(M_3,g_3,\mathcal{F}_3)$ given by $\eta_i\colon (W,g_W)\to (M_i,g_i)$ for $i=2,3$. Now consider the smooth manifold $U\,_{\pi_2}\!\times_{\eta_2}W$ defined as

\begin{equation*}
    U\,_{\pi_2}\!\times_{\eta_2}W:=\left\{(u,w)\in U\times W\,|\, \pi_2(u)=\eta_2(w)\in M_2\right\}
\end{equation*}
\noindent with canonical projections $p_U\colon U\,_{\pi_2}\!\times_{\eta_2}W \to U$ and $p_W\colon U\,_{\pi_2}\!\times_{\eta_2}W \to W$. Note that the tangent space at $(u,w)\in U\,_{\pi_2}\!\times_{\eta_2}W$ is given by

\begin{equation*}
    T_{(u,w)}(U\,_{\pi_2}\!\times_{\eta_2}W)=\left\{(X,Y)\in T_uU\times T_wW\,|\, \exd_u\pi_2(X)=\exd_w\eta_2(Y)\right\}
\end{equation*}
\noindent since every smooth curve on $U\,_{\pi_2}\!\times_{\eta_2}W$ can be expressed as  $(\gamma_U,\gamma_W)$ where $\gamma_U$ and $\gamma_W$ are smooth curves on $U$ and $W$, respectively, such that $\pi_2(\gamma_U)=\eta_2(\gamma_W)$. We now define a Riemannian metric $g$ on $U\,_{\pi_2}\!\times_{\eta_2}W$ as follows:\footnote{We were informed that this idea has been used already in \cite{wink} and \cite{delhoyo}.} 

\begin{equation}
    g((X_1,Y_1),(X_2,Y_2)):=g_U(X_1,X_2)+g_W(Y_1,Y_2)-g_2(\exd_u\pi_2(X_1),\exd_u\pi_2(X_2))
\end{equation}
\noindent where $(X_i,Y_i)\in T_{(u,w)}(U\,_{\pi_2}\!\times_{\eta_2}W)$ for $i=1,2$, and note that $\exd_u\pi_U(X_1)=\exd_w\eta_W(Y_i)$ for $i=1,2$. It is clearly smooth and symmetric. In addition we have
\begin{equation*}
g((X,Y),(X,Y))=\|X\|^2+\|Y\|^2-\|\exd_u\pi_2(X)\|^2=\|X\|^2+\|Y\|^2-\|\exd_w\eta_2(Y)\|^2\geq 0 
\end{equation*}
\noindent for every $(X,Y)\in T_{(u,w)}(U\,_{\pi_2}\!\times_{\eta_2}W)$ since $\pi_U$ and $\eta_W$ are Riemannian submersions, and it is zero if and only if both $X$ and $Y$ are zero vectors. Hence $(U\,_{\pi_2}\!\times_{\eta_2}W,g)$ defines a Riemannian manifold. Now we claim that the projections $p_U$ and $p_W$ are Riemannian submersions. We have
\begin{equation*}
    \ker(\exd_{(u,w)}p_U)=\left\{(0,Y)\in T_uU\times T_wW\,|\, \exd_w\eta_2(Y)=0 \right\} ,
\end{equation*}
 \noindent so its orthogonal complement is given by
\begin{equation*}
    H_{(u,w)}=\left\{(X,Y)\in T_{(u,w)}(U\,_{\pi_2}\!\times_{\eta_2}W)\,|\, g_W(Y,Z)=0\,\,\,\,\,\forall\,Z\in \ker(\exd_w\eta_2) \right\}.
\end{equation*}
\noindent Using the fact that $\eta_W$ is a Riemannian submersion, for every two vectors $(X_1,Y_1)$ and $(X_2,Y_2)$ in $H_{(u,w)}$ we have
\begin{align*}
    g((X_1,Y_1),(X_2,Y_2))&=g_U(X_1,X_2)+g_W(\exd_w\eta_2(Y_1),\exd_w\eta_2(Y_2))-g_2(\exd_w\eta_2(Y_1),\exd_w\eta_2(Y_2))\\ &=g_U(X_1,X_2)=g_U(\exd_{(u,w)}p_U(X_1,Y_1),\exd_{(u,w)}p_U(X_2,Y_2))
\end{align*}
\noindent which proves that $p_U$ is a Riemannian submersion. It has connected fibers since for every $u\in U$, we have $p_U^{-1}(u)=\left\{u\right\}\times\eta_2^{-1}(\pi_2(u))$, which is connected.
Similarly it is shown that $p_W$ is a Riemannian submersion with connected fibers. These two Riemannian submersions are surjective by construction. So the Riemannian manifold $(U\,_{\pi_2}\!\times_{\eta_2}W,g)$ and the surjective Riemannian submersions with connected fibers $\pi_1\circ\p_U$ and $\pi_2\circ\p_W$ define a Morita equivalence between $(N_1,g_1,\mathcal{F}_1)$ and $(N_3,g_3,\mathcal{F}_3)$. This completes the proof.  
\end{proof}

\noindent Although the leaf space of an SRF may not be topologically well-behaved, it inherits a pseudo-metric space structure from the Riemannian metric. Following \cite{PPT10}, for every two leaves $L_1$ and $L_2$ of an SRF $(M,g,\mathcal{F})$, the distance between them is defined by 

\begin{equation*}
    d_{M/\mathcal{F}}(L_1,L_2):=inf\left\{\sum_{i=1}^NL_g(\gamma_i) \right\} .
\end{equation*}
\noindent Here the infimum is taken over all \emph{discrete paths} $(\gamma_1,\cdots,\gamma_N)$ joining $L_1$ and $L_2$, i.e.\ a family of piecewise smooth curves $\gamma_1,\cdots, \gamma_N\colon [0,1]\to M$ for some positive integer $N$, such that $\gamma_1(0)\in L_1$, $\gamma_N(1)\in L_2$ and $\gamma_i(1)$ and $\gamma_{i+1}(0)$ belong to the same leaf for each $i=1,\cdots,n-1$.

\noindent As a corollary of Remark \ref{underlying} and Theorem \ref{hme}, there exists a homeomorphism between the leaf spaces of Morita equivalent module SRFs. The following theorem is the Riemannian counterpart of part $(i)$ of Theorem \ref{hme}:

\begin{theorem}
Let $(N_1,g_1,\mathcal{F}_1)$ and $(N_2,g_2,\mathcal{F}_2)$ be Morita equivalent module SRFs. Then the homeomorphism between the leaf spaces given in Theorem \ref{hme} is distance preserving.
\end{theorem}

\noindent \begin{proof}
Assume that $(M_1,g_1,\mathcal{F}_1)\sim_{ME}(M_2,g_2,\mathcal{F}_2)$ is given by $\pi_i\colon (N,h)\to (M_i,g_i)$ for $i=1,2$. Let $L_1$ and $L_1'$ be two leaves in $(M_1,g_1,\mathcal{F}_1)$ and let $L_2$ and $L_2'$ be their corresponding leaves in $(M_2,g_2,\mathcal{F}_2)$. Consider a discrete path $(\gamma_1,\cdots,\gamma_n)$ joining $L_1$ and $L_1'$. By lifting each $\gamma_i$ into finitely many piecewise smooth horizontal paths, one obtains a discrete path $(\eta_1,\ldots,\eta_{n'})$ for some $n' \geq n$ on $n$ joining $\pi_1^{-1}(L_1)$ and $\pi_1^{-1}(L_1')$ with the same length as $(\gamma_1,\cdots,\gamma_n)$---since the lifts are horizontal with respect to the Riemannian submersion  $\pi_1$. Since $\pi_2$ is a Riemannian submersion, $(\pi_2(\eta_1),\cdots,\pi_2(\eta_{n'}))$ is a discrete path joining $L_2$ and $L_2'$ with a length which is smaller than or equal to the length of $(\gamma_1,\cdots,\gamma_n)$---since the lifts are not necessarily horizontal with respect to $\pi_2$. Consequently \[d_{M_1/\mathcal{F}_1}(L_1,L_1')\geq d_{M_2/\mathcal{F}_2}(L_2,L_2').\]  Similarly $d_{M_2/\mathcal{F}_2}(L_2,L_2')\geq d_{M_1/\mathcal{F}_1}(L_1,L_1')$, which implies $d_{M_1/\mathcal{F}_1}(L_1,L_1')=d_{M_1/\mathcal{F}_2}(L_2,L_2')$. This proves the statement.
\end{proof}

\noindent To define a category \textbf{SRF} of module SRFs one needs to specify their morphisms. We are not going to do this in the present article. But since any good notion of such morphisms should include Riemannian submersions which satisfy that the pullback of the SF on the base agrees with the SF on the total space, we define the following full subcategory $\mathbf{SRF}_0$:

\begin{Definition}
The category $\mathbf{SRF}_0$ has module SRFs as its objects and Riemannian submersions $\pi\colon(N,h,\mathcal{F}_N)\to(M,g,\mathcal{F}_M)$  satisfying $\pi^{-1}\mathcal{F}_M=\mathcal{F}_N$ as its morphisms.
\end{Definition}

\section{\texorpdfstring{$\cI$}{I}-Poisson manifolds} \label{i-poisson}

\noindent In what follows $(P,\left\{\cdot,\cdot\right\})$ stands for a manifold $P$ equipped with an $\R$-bilinear bracket $\left\{\cdot,\cdot\right\} \colon \bigwedge^2 C^\infty(P) \to  C^\infty(P)$ satisfying the Leibniz rule. In other words, $\left\{f,g\right\}= \Pi(\exd f, \exd g)$ for some bivector field $\Pi \in \Gamma(\bigwedge^2TP)$. 
The bracket does not necessarily satisfy the Jacobi identity; if it does, $(P,\left\{\cdot,\cdot\right\})$ is a Poisson manifold \cite{dasilva,LGPV13,CFM21}.  By abuse of notation, we denote the restriction of the bracket to any open subset $U\subset P$ simply by $\left\{\cdot,\cdot\right\}$. Given a function $H\in C^{\infty}(P)$, we call $X_H := \{ H, \cdot \}$ 
the Hamiltonian vector field of $H$ and denote its  flow by $\Phi_H^t$. 
We denote the sheaf of smooth functions on $P$ by $C^\infty$.

\begin{Definition}
A subsheaf $\mathcal{I}$ of a sheaf of rings $\mathcal{O}$ on a manifold $P$
is called \emph{locally finitely generated} if for every $q\in P$ there exist an open neighborhood $U\subset P$ containing $q$ and finitely many sections 
$g_1,...,g_N \in \mathcal{O}(U)$ 
such that $\mathcal{I}(V)=\la g_1\vert_V,...,g_N\vert_V\ra_{\mathcal{O}(V)}$ for every open subset $V\subset U$. 
\end{Definition}

\begin{Remark}
In this article we mostly work with $\mathcal{O}$ being the sheaf of smooth functions, but one may equally consider sheaves of polynomial, real analytic or holomorphic functions for the appropriate choice of $P$.
\end{Remark}

\begin{Definition}\label{ipoisson}
An \emph{$\mathcal{I}$-Poisson manifold} is a triple $(P,\left\{\cdot,\cdot\right\},\mathcal{I})$ where $\mathcal{I}$ is a locally finitely generated subsheaf of smooth functions on $P$, such that for every open subset $U\subset P$ we have
\begin{enumerate}
    \item $\mathcal{I}(U)$ is a $C^\infty(U)$-module,\label{ipcond1}
    \item $\mathcal{I}(U)$ is closed under the bracket,\label{ipcond2}
    \item $ \{\{f,g\},h\}+\{\{g,h\},f\}+\{\{h,f\},g\}\in \mathcal{I}(P),\qquad \forall\, f,g,h\in N(\mathcal{I}),$\label{ipcond3}
\end{enumerate}
where $N(\mathcal{I}):=\left\{f\in C^{\infty}(P)\colon \left\{f\vert_U,\mathcal{I}(U)\right\}\subset \mathcal{I}(U)\,\,\textit{for every open subset }\, U\right\}$.\vspace{1mm}

\noindent  We call $\left\{\cdot,\cdot\right\}$ the \emph{$\mathcal{I}$-Poisson bracket} and $N(\mathcal{I})$ the \emph{$\mathcal{I}$-Poisson normalizer}. 
\end{Definition}

\begin{Example}
Every Poisson manifold $(P,\left\{\cdot,\cdot\right\})$ is canonically an  \emph{$\mathcal{I}$-Poisson manifold} for $\mathcal{I}$ generated by the zero function. 
\end{Example}

\begin{Example} \label{simple}
Let  $(P,\left\{\cdot,\cdot\right\})$ be a Poisson manifold and $I = \langle f_1, \ldots , f_N\rangle_{C^{\infty}(P)}$ be a finitely generated ideal of $C^{\infty}(P)$ which is a Poisson subalgebra. Then the sheaf $\mathcal{I}$  defined by 
\[ U \mapsto \mathcal{I}(U):= \langle f_1\vert_U , \ldots , f_N\vert_U \rangle_{C^{\infty}(U)} \] 
defines an $\mathcal{I}$-Poisson manifold.
\end{Example}

\begin{Example}\label{cored}
Let $(P,\left\{\cdot,\cdot\right\})$ be a Poisson manifold and $C\subset P$ an embedded coisotropic submanifold. Then the triple $(P,\left\{\cdot,\cdot\right\},\mathcal{I}_C)$ where $\mathcal{I}_C(U):=\left\{f\in C^{\infty}(U)\,\colon\,f|_{C\cap U}\equiv0\right\}$ for every open subset $U\subset P$ defines an $\mathcal{I}$-Poisson manifold. Note that in this example, $\mathcal{I}_C$
is in general not finitely generated, only locally so. 
\end{Example}

\begin{Example} \label{quasi-poisson}
Let $(P,\{\cdot,\cdot\},\Phi)$ be a Hamiltonian quasi-Poisson manifold \cite{anton}: $P$ is a $G$-manifold for a compact Lie group $G$, $\left\{\cdot,\cdot\right\} \colon \bigwedge^2 C^\infty(P) \to  C^\infty(P)$ is an $\R$-bilinear bracket  satisfying the Leibniz rule, such that

\begin{equation}\label{almjac}
        \{\{f,g\},h\}+\{\{g,h\},f\}+\{\{h,f\},g\}=\phi_P(df,dg,dh)\qquad \forall f,g,h\in C^\infty(P)
\end{equation}
where $\phi_P\in \mathfrak{X}^3(P)$ is a $3$-vector field induced by the cartan $3$-tensor $\phi\in \bigwedge^3\mathfrak{g}$, and $\Phi \colon P \to G$ a $G$-equivariant map satisfying the \emph{moment map condition}

\begin{equation}\label{mommap}
    \{\Phi^*f,\cdot\}=\tfrac{1}{2}\Phi^*((e_a^L+e_a^R)\cdot f) (e_a)_P\,\qquad \forall f\in C^\infty(G),
\end{equation}
where $(e_a)$ is a basis for $\mathfrak{g}$, $e_a^L,e_a^R\in \mathfrak{X}(G)$ are the left-invariant and and right-invariant vector fields associated to $e_a$ respectively, and $(e_a)_P\in \mathfrak{X}(P)$ is the fundamental vector field induced by $e_a$. 
    
\noindent Let $P_*$ be the open subset of $P$ on which $G$ acts freely.
Fix a conjugacy class $C\subset G$ and let $\mathcal{I}\subset C^\infty(P)$ to be defined as the vanishing ideal of $C_*:=\Phi^{-1}(C)\cap P_*$. Now the triple $(P_*,\{\cdot,\cdot\},\mathcal{I})$ defines an $\mathcal{I}$-Poisson manifold: In Definition \ref{ipoisson}, Condition \ref{ipcond1} is clear and Condition \ref{ipcond2} is a consequence of the moment map condition and $G$-equivariance of $\Phi$. It remains to show that Condition \ref{ipcond3} is satisfied: The moment map condition gives

\begin{equation*}
        \mathcal{N}(\mathcal{I})=\{f\in C^\infty(P)\,\colon\, f\vert_{C_*}\in C^\infty(C_*)^G\},
\end{equation*}
which implies that for $f,g,h\in \mathcal{N}(\mathcal{I})$, the function $\phi_P(df,dg,dh)$ vanishes on $C_*$ since $\phi_P\vert_{C_*}\in \mathfrak{X}^3(C_*)$.
\end{Example}

\noindent Under some conditions, a reduction process applied to Examples \ref{cored} and \ref{quasi-poisson} results in \emph{reduced Poisson manifolds}.

\begin{Example}\label{cored2}
    In Example \ref{cored}, the Hamiltonian vector fields of functions in $\mathcal{I}_C$ are tangent to $C$ and they are closed under the Lie bracket, hence defining an SF on $C$. If this SF is regular and the quotient map $\pi\colon C\to C_{red}$ to the leaf space $C_{red}$ is a smooth submersion, then $C_{red}$ inherits a Poisson bracket $\left\{\cdot,\cdot\right\}_{red}$ such that $\pi^*\left\{f,g\right\}_{red}=\left\{F,G\right\}|_C$, where $F$ and $G$ are smooth functions on $P$ satisfying $F|_C=\pi^*f$ and $G|_C=\pi^*g$. This process is called the \emph{coisotropic reduction} \cite{marsdenratiu}.
\end{Example}

\begin{Example}\label{quasi-poisson2}
    In Example \ref{quasi-poisson}, Theorem $6.1$ in \cite{anton} implies that the quotient $C_{red}:=C_*/G$ inherits a Poisson bracket $\{\cdot,\cdot\}_{red}$.
\end{Example}

\begin{Remark}
The notion of an $\mathcal{I}$-Poisson manifold is motivated by generalizing Examples \ref{cored} and \ref{quasi-poisson} and their reductions to a potentially singular setting, where the quotient $C_{red}$ does not need to exist as a manifold and the reduction is performed algebraically.
\end{Remark}

\noindent As a consequence of Conditions \ref{ipcond2} and \ref{ipcond3} of Definition \ref{i-poisson}, the quotient $N(\mathcal{I})/\mathcal{I}(P)$ forms a Poisson algebra. This motivates the following definition:

\begin{Definition}\label{reducedPois}
The \emph{reduced Poisson algebra} of the $\mathcal{I}$-Poisson manifold $(P,\left\{\cdot,\cdot\right\},\mathcal{I})$ is defined to be the Poisson algebra $\mathcal{R}(\mathcal{I}):=N(\mathcal{I})/\mathcal{I}(P)$.
\end{Definition}

\begin{Remark} This is a straightforward generalization of the set of Dirac observables \cite{Dirac}. 
The algebra $\mathcal{R}(\mathcal{I})$ also appears  in \cite{sw83} as an algebraic method of reducing Hamiltonian $G$-spaces with singular moment maps.  
\end{Remark}

\begin{Example}
 If $C_{red}$ in Example \ref{cored2}  is smooth, then $\mathcal{R}(\mathcal{I}_C)$ is isomorphic to the Poisson algebra $C^{\infty}(C_{red})$. Similarly, in Example \ref{quasi-poisson2}, the Poisson algebra of functions on $C_{red}$ is isomorphic to the Poisson algebra $\mathcal{R}(\mathcal{I})$.
\end{Example}

\begin{Example}\label{action}
Let $G$ be a connected Lie group acting on a Poisson manifold $(P,\left\{\cdot,\cdot\right\})$ by Poisson diffeomorphisms with a $G$-equivariant moment $\mu\colon P\rightarrow \mathfrak{g}^*$. Following \cite{sw83}, the subsheaf $\mathcal{I}\subset C^\infty$ generated by smooth functions $\langle \mu,\mathfrak{g}\ra$ is a Poisson subalgebra and one has \[\mathcal{R}(\mathcal{I})\cong (C^\infty(P)/\mathcal{I})^G.\] Moreover, if $G$ is compact, then Proposition $5.12$ in \cite{gotay} states that\[\mathcal{R}(\mathcal{I})\cong C^\infty(P)^G/\mathcal{I}^G.\] 
\end{Example}

\begin{Example}
Let $P=T^*\R^n$, $n>1$, with coordinates $(q^1,\ldots,q^n,p_1,\ldots,p_n)$ and $\mathcal{I}\subset C^\infty$ the subsheaf generated by the $n(n-1)/2$ functions $q^ip_j-q^jp_i$ for $1\leq i< j\leq n$. This is a special case of Example \ref{action} for the diagonal action of $G=SO(n)$ on $T^*\R^n$. We have \[\mathcal{R}(\mathcal{I})\cong W^\infty(D),\] where $D\subset \R^3$ is defined by 
\[D:=\left\{(x_1,x_2,x_3)\in \R^3\,\vert\,x_1^2+x_2^2=x_3^2\,\,\,\textit{and}\,\,\,x_3\geq0\right\}\] and $W^\infty(D)$ stands for the smooth functions on $D$ in the sense of Whitney, \emph{i.e.}\ the restriction of $C^\infty(\R^3)$ to $D$. For more details and proofs see Theorem $5.6$ and Example $5.11(a)$ of \cite{gotay}. 

\noindent  The Poisson bracket on $W^\infty(D)$ can be understood as follows: Identify $\R^3$ with the Poisson manifold $\mathfrak{so}(2,1)^*$ and, simultaneously, with 2+1 dimensional Minkowski space. The symplectic leaves of $\mathfrak{so}(2,1)^*$ then consist of spacelike vectors of a fixed Minkowski norm (one-sheeted hyperboloids), null vectors decompose into the origin, the forward light cone, and the backward light cone as three distinct leaves, and finally timelike vectors of a fixed norm yield two leaves each (two-sheeted hyperboloids). Then restriction to $D$ corresponds precisely to restricting to the forward lightcone and the origin in this Minkowski space. This bracket does not depend on the extension of a function on $D$ to the ambient space since $D$ is the collection of (two) symplectic leaves.
\end{Example}

\begin{Remark}
If in the previous example one restricts to the polynomial functions, such that $\mathcal{I}\subset \R[q^1,\ldots,q^n,p_1,\ldots,p_n]$, one finds 
\[\mathcal{R}(\mathcal{I})\cong S^\bullet \!\left(\mathfrak{so}(2,1)\right) / \la x_1^2 + x_2^2-x_3^2 \ra , \]
\emph{i.e.}\ the polynomial functions on $\mathfrak{so}(2,1)^*$ modulo the ideal generated by the quadratic Casimir. So one looses the restriction $x_3\geq 0$ that one finds in the smooth setting.
\end{Remark}

\begin{Definition}
A \emph{dynamical $\mathcal{I}$-Poisson manifold} denoted by $(P,\left\{\cdot,\cdot\right\},\mathcal{I},H)$ consists of an $\mathcal{I}$-Poisson manifold $(P,\left\{\cdot,\cdot\right\},\mathcal{I})$ and a Hamiltonian function $H\in N(\mathcal{I})$. Its reduction is defined to be the pair $(\mathcal{R}(\mathcal{I}),[H])$ where $[H]\in \mathcal{R}(\mathcal{I})$ is the equivalence class of $H$.
\end{Definition}

\noindent The following proposition reveals one of the main properties of dynamical $\mathcal{I}$-Poisson manifolds.

\begin{Proposition} \label{normcons}
Let $(P,\left\{\cdot,\cdot\right\},\mathcal{I},H)$ be a dynamical $\mathcal{I}$-Poisson manifold. Then the Hamiltonian flow of $H$ locally preserves $\mathcal{I}$, \emph{i.e.}\ for every $q_0\in P$ there exists an open neighborhood $U\subset P$ around $q_0$ such that $\Phi_H^t|_U$ is defined for $t\in(-\epsilon,\epsilon)$ and

\begin{equation}
    ({\Phi_H^t})^* \, \mathcal{I}(\Phi_H^t(U))=\mathcal{I}(U).
\end{equation}
In the case that the Hamiltonian vector field $X_H$ is complete, this implies that, for all $t \in \R$, one has $({\Phi_H^t})^*\,\mathcal{I}\circ \Phi_H^t =\mathcal{I}$  and, in particular, that the ideal $\mathcal{I}(P)$ is preserved,
\begin{equation*}
     ({\Phi_H^t})^* \, \mathcal{I}(P)=\mathcal{I}(P).
\end{equation*}
\end{Proposition}

\noindent\begin{proof}
Choose an open neighborhood $W\subset P$ around $q_0$ where $\mathcal{I}(W)$ is generated by finitely many functions $g_1,\ldots,g_N$ for some positive integer $N$. Then by the existence and uniqueness theorem for ODEs there exist an open subset $U\subset W$ containing $q_0$ and an interval $\left(-\epsilon,\epsilon\right)$, $\epsilon>0$, such that $\Phi_H^t|_U$ is defined for $t\in(-\epsilon,\epsilon)$. By the definition of the $\mathcal{I}$-Poisson normalizer in Definition \ref{ipoisson}, there exist functions $\lambda_a^b\in C^{\infty}(U)$, $a,b=1,\ldots,N$, such that:

\begin{equation*}
\left \{ H,g_a\right \}=\sum_{b=1}^N\lambda_a^b\, g_b.  
\end{equation*}

\noindent Using this equation, 
we obtain:
\begin{equation}\label{inftoint}
\tfrac{d}{dt}( ({\Phi_H^t})^*g_a)= ({\Phi_H^t})^*\left \{H,g_a\right \}= \sum_{b=1}^N (({\Phi_H^t})^*\lambda_a^b)(({\Phi_H^t})^*g_a) .
\end{equation}

\noindent Now, let $X_p(t)\in \R^N$ be a column vector with $a$-th component equal to $g_a\circ \Phi_H^t (p)$ for $a=1,\ldots,N$ and let $A_p(t)$ be the $N$ by $N$ matrix $(\lambda_a^b\circ \Phi_H^t(p))_{a,b=1}^N$. Equation \eqref{inftoint} then transforms into the following family of non-autonomous linear ODEs

\begin{equation}\label{nonautode}
    \tfrac{\exd}{\exd t}X_p(t)=A_p(t)X_p(t).
\end{equation}
This equation and its 
initial conditions  depend smoothly on $p\in U$. It is standard knowledge that solutions  to \eqref{nonautode} take the form:
\begin{equation}\label{solsol}
  X_p(t)=\Psi_p(t)X_p(0) .
\end{equation}
\noindent Here $\Psi_p(t)=(\psi_a^b(t,p))_{a,b=1}^N$ is the fundamental matrix of the ODE, satisfying $\Psi_p(0)=\mathrm{I}_N$ and
\begin{align*}
    \tfrac{\exd}{\exd t}\Psi_p(t)=A_p(t)\Psi_p(t).
\end{align*}
$\Psi_p(t)$ is sometimes also called the (time-) ordered exponential of $A_p(t)$.

\noindent Since $A_p(t)$ and $X_p(0)$ depend smoothly on $p$, the components of the fundamental matrix, $\psi_a^b(t,p)$, depend smoothly on $p$ as well. Now, Equation \eqref{solsol} can be written as 

\begin{equation*}
   ( {\Phi_H^t})^*g_a(p)=\sum_{b=1}^N \psi_a^b(t,p)g_b(p),
\end{equation*}
which implies the inclusion $({\Phi_H^t})^* \, \mathcal{I}(\Phi_H^t(U))\subset\mathcal{I}(U)$. 

\noindent To prove equality, we first observe that the inclusion yields also $({\Phi_H^{-t}})^* \, \mathcal{I}(U)\subset\mathcal{I}(\Phi_H^t(U))$ for every $t\in (-\epsilon,\epsilon)$. Thus, for every $f\in \mathcal{I}(U)$, one has 
$({\Phi_H^{-t}})^*f\in \mathcal{I}(\Phi_H^t(U))$. But on the other hand, 
we have the obvious identity
\[f=({\Phi_H^t})^*\left(({\Phi_H^{-t}})^*f\right),\]
and therefore $f\in ({\Phi_H^t})^* \, \mathcal{I}(\Phi_H^t(U))$.
\end{proof}

\noindent The following example shows that the condition of being locally finitely generated in the definition of $\mathcal{I}$-Poisson manifolds is crucial for Proposition \ref{normcons} to hold true:

\begin{Example}
Consider the Poisson manifold $M=T^*\R\cong \R^2$ with coordinates $(q,p)$ and standard Poisson bracket 

\[\left\{f,g\right\}=\pderiv{f}{p}\pderiv{g}{q}-\pderiv{f}{q}\pderiv{g}{p}.\]

\noindent Let $\mathcal{I}$ be the subsheaf of $C^\infty$ vanishing on $\left\{q<0\right\}\subset M$, which is not locally finitely generated around every point on the $p$-axis, but still closed under the Poisson bracket. Then the coordinate function $p$ is an element of $N(\mathcal{I})$ since $X_p=\pderiv{}{q}$ preserves $\mathcal{I}$. But the Hamiltonian flow of $X_p$ is given by $\Phi_{X_p}^t(q,p)=(q+t,p)$, which evidently does not preserve $\mathcal{I}$ if $t>0$.
\end{Example}

\noindent In order to define the category of $\mathcal{I}$-Poisson manifolds, we introduce a notion of morphisms and show that they can be composed:

\begin{Definition}\label{ipoissoncat}
Let $\varphi\colon(P_1,\left\{\cdot,\cdot\right\}_1,\mathcal{I}_1)\to(P_2,\left\{\cdot,\cdot\right\}_2,\mathcal{I}_2)$ be a smooth map between two $\mathcal{I}$-Poisson manifolds. We call it an $\mathcal{I}$-Poisson map if the following three conditions are satisfied:
\begin{align}
   \varphi^*(\mathcal{I}_2(P_2))&\subset& \mathcal{I}_1(P_1)\label{first}, \\
   \varphi^*N(\mathcal{I}_2)&\subset& N(\mathcal{I}_1)\label{dirobs},\\
   \{\varphi^* f,\varphi^* g\}_1 - \varphi^* \{f,g\}_2 &\in&  \mathcal{I}_1(P_1) \qquad \forall f,g \in N(\mathcal{I}_2)\label{third}.
\end{align}
For dynamical $\mathcal{I}$-Poisson manifolds we add the condition $\varphi^* H_2-H_1\in  \mathcal{I}_1$.
\end{Definition}

\begin{Proposition}
The composition of two $\mathcal{I}$-Poisson maps is an $\mathcal{I}$-Poisson map.
\end{Proposition}

\noindent \begin{proof}
Consider the following $\mathcal{I}$-Poisson maps:

\begin{align*}
    \varphi\colon(P_1,\left\{\cdot,\cdot\right\}_1,\mathcal{I}_1)&\to&(P_2,\left\{\cdot,\cdot\right\}_2,\mathcal{I}_2)\\
    \psi\colon(P_2,\left\{\cdot,\cdot\right\}_2,\mathcal{I}_2)&\to&(P_3,\left\{\cdot,\cdot\right\}_3,\mathcal{I}_3).
\end{align*}
Equations \eqref{first} and \eqref{dirobs} for 
$\psi\circ\varphi$ follow directly from those equations for $\psi$ and $\varphi$. 
It is thus enough to verify Equation \eqref{third} for the composition. For all $f,g \in N(\mathcal{I}_3)$ we have
\begin{align*}
    &\quad\left\{f\circ \psi\circ\varphi,g\circ\psi\circ \varphi\right\}_1-\left\{f,g\right\}_3\circ\psi\circ \varphi\\
    &=\left\{\left(f\circ \psi\right)\circ\varphi,\left(g\circ\psi\right)\circ \varphi\right\}_1-\left\{f\circ \psi,g\circ\psi\right\}_2\circ\varphi\\
    &+&\left(\left\{f\circ \psi,g\circ\psi\right\}_2-\left\{f,g\right\}_3\circ\psi\right)\circ\varphi\\
    &\in& \mathcal{I}_1(P_1)+\varphi^*\mathcal{I}_2(P_2)\subset \mathcal{I}_1(P_1),
\end{align*}
where we used Equations \eqref{first} and \eqref{third} for $\varphi$ and Equations \eqref{dirobs} and \eqref{third} for $\psi$ in the last line of the proof. A similar computation shows that morphisms of dynamical $\mathcal{I}$-Poisson manifolds can be composed as well.
\end{proof}

\begin{Definition}\label{semistrict}
The category \emph{\textbf{IPois}} and \emph{\textbf{dynIPois}} consist of $\mathcal{I}$-Poisson manifolds together with $\mathcal{I}$-Poisson maps and dynamical $\mathcal{I}$-Poisson manifolds together with dynamical $\mathcal{I}$-Poisson maps, respectively. By requiring the $\mathcal{I}$-Poisson bracket to be a Poisson bracket, we obtain a subcategory which we call \emph{(dynamical) semi-strict $\mathcal{I}$-Poisson manifolds} \emph{\textbf{ssIPois}} \emph{(\textbf{ssdynIPois})}. Similarly, the category \emph{\textbf{sIPois}} \emph{(\textbf{sdynIPois})} of strict \emph{(dynamical)  $\mathcal{I}$-Poisson manifolds}  is defined by requiring that the $\mathcal{I}$-Poisson bracket is a Poisson bracket and that the morphisms are  Poisson maps. 
\end{Definition}

\begin{Remark} While \textbf{ssIPois} is a full subcategory of \textbf{IPois}, \textbf{sIPois} is not.
\end{Remark}

\begin{Remark} \label{minimal}
The three conditions in Definition \ref{ipoissoncat} are the minimal conditions for the map $\varphi^*$ to induce a morphism of Poisson algebras $\Tilde{\varphi}\colon \mathcal{R}(\mathcal{I}_2)\to \mathcal{R}(\mathcal{I}_1)$. In particular, we obtain a functor $F$ from \textbf{IPois} $\!{}^{op}$ to  \textbf{PoisAlg}, the category of Poisson algebras. We call $F$ the reduction functor. 
\end{Remark} 

\begin{Remark} \label{subcategory}
Viewing Poisson manifolds $(P, \left\{\cdot,\cdot\right\})$ as   $\mathcal{I}$-Poisson manifolds $(P, \left\{\cdot,\cdot\right\},0)$,  $\mathcal{I}$-Poisson maps are precisely Poisson maps. This identifies the category of Poisson manifolds \textbf{Pois} with a full subcategory of \textbf{sIPois}.
\end{Remark}

\begin{Remark} \label{Waldmann}
There is a functor from  $\mathbf{IPois}$ to $\mathbf{C_3Alg}$, the category of \emph{coisotropic triples of algebras} as introduced in \cite{DEW19}. On the level of objects, one assoicates the triple $(C^\infty(P),N(\mathcal{I}),\mathcal{I}(P))$ to 
 every $\mathcal{I}$-Poisson manifold $(P,\left\{\cdot,\cdot\right\},\mathcal{I})$, while a  morphism $\varphi$ in our sense gives rise to a morphism $\varphi^*$ in $\mathbf{C_3Alg}$ due to the first two defining conditions  \eqref{first} and \eqref{dirobs}.
\end{Remark}

\section{Singular (Riemannian) foliations through \texorpdfstring{$\cI$}{I}-Poisson manifolds}\label{objects}

\noindent Let $M$ be a smooth manifold. We denote by $C_k^{\infty}(T^*M)\subset C^{\infty}(T^*M)$  the algebra of homogeneous polynomials of degree $k$ in the fiber coordinates of $T^*M$ with coefficients in $C^{\infty}(M)$. Every vector field $X\in \mathfrak{X}(M)$ defines an element $\overline{X}\in C_1^{\infty}(T^*M)$ on the cotangent bundle of $M$, defined as

\begin{equation*}
\overline{X}(q,p)=\langle p, X\vert_q \rangle
\end{equation*}
\noindent for every $(q,p)\in T^*M$, where $q\in M$ and $p\in T_q^*M$ and $\langle \cdot , \cdot \rangle$ denotes the canonical pairing. This construction can be naturally extended to the sections of $\Gamma(S^k(TM))$ to obtain elements in $C_k(T^*M)$.

\begin{Lemma} \label{liepoisson}
Let $X,Y\in \mathfrak{X}(M)$ be two vector fields on $M$. Then

\begin{equation}
\left \{ \overline{X},\overline{Y}\right\}_{T^*M}=\overline{[X,Y]},
\end{equation}
 where $\left\{\cdot,\cdot\right \}_{T^*M}$ is the canonical Poisson bracket on $T^*M$ and $[\cdot,\cdot]$ is the Lie bracket of vector fields on $M$.
\end{Lemma}

\noindent \begin{proof}
Let $(q^1,...,q^n)$ be a local coordinate system on $M$, and $(q^1,...,q^n,p_1,...,p_n)$ the corresponding canonical local coordinates on $T^*M$. In this coordinate system $\overline{X}$ and $\overline{Y}$ can be written as $\overline{X}=\sum_{i=1}^n X^ip_i$ and $\overline{Y}=\sum_{i=1}^n Y^i p_i$, where  $X^i,Y^i$ are the components of $V,W$ in the above coordinate system. The following calculation proves the lemma:

\begin{equation*}
\left\{\overline{X},\overline{Y}\right \}_{T^*M}=\sum_{i=1}^n(\sum_{j=1}^nX^j\tfrac{\partial Y^i}{\partial q^j}-Y^j\tfrac{\partial X^i}{\partial q^j})p_i=\sum_{i=1}^n[X,Y]^ip_i=\overline{[X,Y]}. 
\end{equation*}
\end{proof}

\noindent Lemma \ref{liepoisson} and the Leibniz rule for the Lie derivative of tensor fields imply:

\begin{Corollary} \label{symvect}
Let $S$ be an element of $\Gamma(S^k(TM))$ for some $k\geq 0$, and $\overline{S}$ be its corresponding element in $C_k^{\infty}(T^*M)$. Then for every vector field $X\in \mathfrak{X}(M)$ we have:
\begin{equation}
\left\{\overline{X},\overline{S}\right \}_{T^*M}=\overline{\mathcal{L}_XS}  .  
\end{equation}
\end{Corollary}

\noindent Now let $(M,\mathcal{F})$ be a foliated manifold. Define a $C^\infty(T^*M)$-submodule $\mathcal{J}_{\mathcal{F}}\subset C^\infty_c(T^*M)$ by
\[\mathcal{J}_{\mathcal{F}}:=\la \overline{X}\,\colon\, X\in \mathcal{F}\ra_{C^\infty_c(T^*M)}.\]
Note that the generators of $\mathcal{J}_{\mathcal{F}}$ are not required to be compactly supported on $T^*M$.

\noindent Then we define the sub-presheaf $\mathcal{I}_{\mathcal{F}}$ of the sheaf of smooth functions on $T^*M$ by
\begin{align}\label{IF}
    \mathcal{I}_{\mathcal{F}}(U):=\left\{f\in C^\infty(U)\,\colon\,\rho f\in \mathcal{J}_{\mathcal{F}}\quad\forall\,\rho\in C^\infty_c(U)\right \}\,
\end{align}
for every open subset $U\subset T^*M$.

\begin{Proposition}
The presehaf $\mathcal{I}_{\mathcal{F}}$ defined in Equation \eqref{IF} is a subsheaf of the sheaf of smooth functions on $T^*M$.
\end{Proposition}

\noindent \begin{proof}
The locality of $\mathcal{I}_{\mathcal{F}}$ is evident, since $\mathcal{I}_{\mathcal{F}}$ is a sub-presheaf of the sheaf of smooth functions on $T^*M$. To verify the gluing property, let $\left\{U_i\right \}_{i=1}^\infty$ be an arbitrary open cover of $T^*M$ and let $f\in C^\infty(T^*M)$ be such that $f\vert_{U_i}\in \mathcal{I}_{\mathcal{F}}(U_i)$ for every positive integer $i$. We prove that $f\in \mathcal{I}_{\mathcal{F}}(T^*M)$ as follows: it is enough to show that for every $\rho\in C^\infty_c(T^*M)$, we have $\rho f\in \mathcal{J}_{\mathcal{F}}$. Since $\mathrm{supp}(\rho)$ is compact, it can be covered by finitely many open subsets $U_{i_1},\ldots,U_{i_N}$ in $\left\{U_i\right \}_{i=1}^\infty$. Choose a partition of unity $\sigma_0,\sigma_{i_1},\ldots,\sigma_{i_N}$ subordinate to the open cover $\left\{U_0:=T^*M\setminus \mathrm{supp}(\rho),U_{i_1},\ldots,U_{i_N}\right \}$ of $T^*M$ and write
\[\rho f=\sum_{k=1}^N\rho\sigma_{i_k}f\vert_{U_{i_k}}.\]
The latter implies that $\rho f\in \mathcal{J}_{\mathcal{F}}$, since by definition of $\mathcal{I}_{\mathcal{F}}(U_{i_k})$, for each $k=1,\ldots,N$ we have $\rho\sigma_{i_k}f\vert_{U_{i_k}}\in \mathcal{J}_{\mathcal{F}}$.
\end{proof}

\noindent We prove that the sheaf $\mathcal{I}_{\mathcal{F}}$ satisfies the properties of Definition \ref{ipoisson}, in the following lemmas:

\begin{Lemma}
For every open subset $U\subset T^*M$ we have 
\begin{equation*}
    \left \{ \mathcal{I}_{\mathcal{F}}(U),\mathcal{I}_{\mathcal{F}}(U)\right \}_{T^*M}\subset \mathcal{I}_{\mathcal{F}}(U)\,
\end{equation*}
\end{Lemma}

\noindent \begin{proof}
Let $f,g\in\mathcal{I}_{\mathcal{F}}(U)$. It is enough to show that for every $\rho\in C^\infty_c(U)$ we have $ \rho\left \{ f,g\right \}\subset \mathcal{J}_{\mathcal{F}}$. Choose a compactly supported function $\sigma \in C^\infty_c(U)$ such that $\sigma\vert_{\mathrm{supp}(\rho)}\equiv 1$. One obtains
\begin{equation*}
    \rho\left\{f,g\right \}_{T^*M} =
    \left \{\sigma f,\rho g\right \}_{T^*M}-\left \{\sigma f,\rho\right \}_{T^*M}g-\rho f\left\{\sigma,g\right \}_{T^*M}\in \mathcal{J}_{\mathcal{F}},
\end{equation*}
since the first term belongs to $\mathcal{J}_{\mathcal{F}}$ by Lemma \ref{liepoisson}, the second term is inside $\mathcal{J}_{\mathcal{F}}$ by Definition of $\mathcal{I}_{\mathcal{F}}(U)$, and the last term vanishes identically. Consequently $\left\{f,g\right \}_{T^*M}\in\mathcal{I}_{\mathcal{F}}(U)$.

\end{proof}

\begin{Lemma}\label{lfgif}
Let $U\subset M$ be an open subset such that $\iota_U^{-1}\mathcal{F}=\la X_1,\ldots,X_N\ra_{C_c^\infty(U)}$ for finitely many vector fields $X_1,\ldots,X_N\in \mathfrak{X}(U)$. Then 
\begin{equation*}
    \mathcal{I}_{\mathcal{F}}(V)=\la\overline{X_1}\vert_V,\ldots,\overline{X_N}\vert_V\ra_{C^\infty(V)},
\end{equation*}
for every open subset $V\subset T^*U$.
\end{Lemma}

\noindent \begin{proof}
We first prove that $\la\overline{X_1}\vert_V,\ldots,\overline{X_N}\vert_V\ra_{C^\infty(V)}\subset \mathcal{I}_{\mathcal{F}}(V)$. Let $\sum_{a=1}^N\lambda^a\overline{X_a}\vert_V$ be an element of $\la\overline{X_1}\vert_V,\ldots,\overline{X_N}\vert_V\ra_{C^\infty(V)}$ and take an arbitrary $\rho\in C^\infty_c(V)$. By choosing a compactly supported function $h\in C^\infty_c(U)$ such that $h\vert_{\mathrm{supp}(\rho)}\equiv 1$ (when viewing $h$ as an element of $C^\infty_0(T^*U)$), we have
\begin{equation*}
    \rho \sum_{a=1}^N\lambda^a\overline{X_a}\vert_V=\sum_{a=1}^N\rho\lambda^a\overline{hX_a}\in \mathcal{J}_{\mathcal{F}},
\end{equation*}
since $\rho\lambda^a\in C_c^\infty(V)$ and $hX_a\in \mathcal{F}$ for all $a=1,\ldots,N$. To prove equality, let $f\in \mathcal{I}_{\mathcal{F}}(V)$. Choose a partition of unity $\left\{\rho_i\right \}_{i=1}^\infty$ subordinate to a locally finite cover $\left\{V_i\right \}_{i=1}^\infty$ of $V$. Since for every $i$ we have $\rho_if\in \mathcal{J}_{\mathcal{F}}$ and $V\subset T^*U$, there exist functions $\lambda_i^1,\ldots,\lambda_i^N\in C^\infty_c(T^*U)$ such that
\begin{equation*}
    \rho_if=\sum_{a=1}^N\lambda_i^a\vert_{V}\overline{X_a}\vert_V.
\end{equation*}
This implies
\begin{align*}
    f=\sum_{i=1}^\infty\rho_if&=\sum_{i=1}^\infty\sum_{a=1}^N\lambda_i^a\vert_V\overline{X_a}\vert_V\\
    &=\sum_{a=1}^N\left(\sum_{i=1}^\infty\lambda_i^a\vert_V\right)\overline{X_a}\vert_V,
\end{align*}
which completes the proof.
\end{proof}

\begin{Corollary}
Let $(M,\mathcal{F})$ be a foliated manifold. Then the triple $(T^*M,\left\{\cdot,\cdot\right \}_{T^*M},\mathcal{I}_{\mathcal{F}})$ is an $\mathcal{I}$-Poisson manifold.
\end{Corollary}

\noindent For every Riemannian manifold $(M,g)$, its cotangent bundle $T^*M$ carries a natural Hamiltonian function $H_g$: 

\begin{equation*}
H_g(q,p)= \tfrac{1}{2} \langle p , g_{\flat}^{-1}(p) \rangle 
\end{equation*}
\noindent for every $(q,p) \in T^*M$, where $g_{\flat} \colon T_qM \to T_q^*M$ is the musical isomorphism $v\mapsto g(v,\cdot)$. In local Darboux coordinates this becomes 
$H_g(q^1,...,q^n,p_1,...,p_n)=\tfrac{1}{2}\sum_{i,j=1}^ng^{ij}p_ip_j$ where the matrix $(g^{ij})$ is the inverse to the matrix of the Riemannian metric $(g_{ij})$ in the coordinates $(q^1,...,q^n)$. Equivalently, we can define also $H_g$ using the isomorphism $\Gamma(S^2(TM))\cong C_2^{\infty}(T^*M)$, under which it becomes identified with $g^{-1}:= \sum_{i,j=1}^n g^{ij} \,  \partial_i \odot \partial_j$, i.e.\ $H_g = \tfrac{1}{2}\overline{g^{-1}}$.

\noindent The following fact about $H_g$ is standard knowledge, which we still prove for completeness.  

\begin{Proposition}\label{Hamflow}The Hamiltonian flow of $H_g$ is the image of the geodesic flow under the musical isomorphism, \emph{i.e.}\ for every geodesic $\gamma\colon(-\epsilon,\epsilon)\to M$ and every $t\in (-\epsilon,\epsilon)$, we have
\[\Phi_{H_g}^t(\gamma(0),g(\dot{\gamma}(0),\cdot))=(\gamma(t),g(\dot{\gamma}(t),\cdot)).\]
\end{Proposition}

\noindent \begin{proof}
Assume that $(q^1,\ldots,q^n)$ is a normal coordinate system centered at $q\in M$, i.e.\ $g_{ij}(q)=\delta_{ij}$ and $\partial_kg_{ij}(q)=0$ for  $i,j,k=1,\ldots,n$. For every $p\in T_q^*M$ we have
\[X_{H_g}(q,p)=\sum_{i=1}^np_i\pderiv{}{q^i}|_q.\]
Let $\gamma\colon(-\epsilon,\epsilon)\to M$ be a geodesic passing through $q$ at $t=0$; in particular, $\ddot{q}^i(0)=0$. Then $((\gamma(t),g(\dot{\gamma}(t),\cdot))$ is a curve on $T^*M$ passing through $(q,p)=((\gamma(0),g(\dot{\gamma}(0)),\cdot)$ at $t=0$; in local coordinates, $((\gamma(t),g(\dot{\gamma}(t)),\cdot)=(q^1(t),\ldots,q^n(t),p_1(t),\ldots,p_n(t))$ where $p_i(t)=\sum_{j=1}^ng_{ij}(q(t))\dot{q}^j(t)$.  Then, since $\dot{p}_i(0)=0$, we have
\[\tfrac{\exd}{\exd t}|_{t=0}((\gamma(t),g(\dot{\gamma}(t),\cdot))=\sum_{i=1}^n\dot{q}^i(0)\pderiv{}{q^i}|_q.\]
On the other hand, $\dot{q}^i(0)=g(\dot{\gamma}(0),\pderiv{}{q^i}|_q)=p_i(0)$, which indeed gives
\[X_{H_g}(q,p)=\tfrac{\exd}{\exd t}|_{t=0}((\gamma(t),g(\dot{\gamma}(t),\cdot)).\]
\end{proof}

\begin{Lemma}\label{normif}
Let $(M,\mathcal{F})$ be a foliated manifold. We have

\begin{equation*}
    N(\mathcal{I}_{\mathcal{F}})=\left\{f\in C^\infty(T^*M)\,\colon\,\left\{f,\mathcal{J}_{\mathcal{F}}\right \}_{T^*M}\subset \mathcal{J}_{\mathcal{F}}\right \}.
\end{equation*}
\end{Lemma}

\noindent\begin{proof}
The inclusion $N(\mathcal{I}_{\mathcal{F}})\subset \left\{f\in C^\infty(T^*M)\,\colon\,\left\{f,\mathcal{J}_{\mathcal{F}}\right \}_{T^*M}\subset \mathcal{J}_{\mathcal{F}}\right \}$ is satisfied by Definition \ref{ipoisson} and the fact that $\mathcal{J}_{\mathcal{F}}$ is equal to the set of compactly supported elements in $\mathcal{I}_{\mathcal{F}}(T^*M)$. Now let $f\in C^\infty(T^*M)$ be such that $\left\{f,\mathcal{J}_{\mathcal{F}}\right \}_{T^*M}\subset \mathcal{J}_{\mathcal{F}}$. Let $U\subset T^*M$ be an open subset and $g\in \mathcal{I}_{\mathcal{F}}(U)$. For every $\rho\in C^\infty_c(U)$ we have

\begin{align*}
    \rho \left\{f\vert_U,g\right \}_{T^*M}=\left\{f,\rho g\right \}_{T^*M}-\left\{f,\rho\right \}_{T^*M}g\in \mathcal{J}_{\mathcal{F}},
\end{align*}
since $\rho g\in \mathcal{J}_{\mathcal{F}}$ and $\left\{f,\rho\right \}_{T^*M}$ is compactly supported in $U$. The definition of $\mathcal{I}_{\mathcal{F}}(U)$ then implies that $\left\{f\vert_U,\mathcal{I}_{\mathcal{F}}(U)\right \}_{T^*M}\subset \mathcal{I}_{\mathcal{F}}(U)$. Since $U$ is arbitrary, we obtain $f\in N(\mathcal{I}_{\mathcal{F}})$.
\end{proof}

\noindent Now we can state an equivalent definition of module SRFs through $\mathcal{I}$-Poisson geometry.

\begin{Proposition} \label{srfcons}
A singular foliation $\mathcal{F}$ on a Riemannian manifold $(M,g)$ defines a module SRF, if and only if

\begin{equation}
H_g\in N(\mathcal{I}_{\mathcal{F}}).
\end{equation}
\end{Proposition}

\noindent \begin{proof}
Assume that $(M,g,\mathcal{F})$ is a module SRF. By Lemma \ref{lem1}, for every $X\in \mathcal{F}$ we have
\[\mathcal{L}_Xg^{-1}\in \mathfrak{X}(M)\odot \mathcal{F}.\]
Using the isomorphism $\Gamma(S^2(TM))\cong C_2^{\infty}(T^*M)$ and Corollary \ref{symvect}, we obtain 

\begin{equation*}
  \left \{
\overline{X},H_g\right \}_{T^*M}=\tfrac{1}{2}\overline{\mathcal{L}_Xg^{-1}}\, \in\,\overline{\mathfrak{X}(M)\,\odot\, \mathcal{F}},
\end{equation*}
which together with the Leibniz rule imply $\left \{
\mathcal{J}_{\mathcal{F}},H_g\right \}_{T^*M}\subset \mathcal{J}_{\mathcal{F}}$. Lemma \ref{normif} then implies that $H_g\in N(\mathcal{I}_{\mathcal{F}})$. Conversely assume that $H_g\in N(\mathcal{I}_{\mathcal{F}})$. After Proposition \ref{locality} we can assume that $\mathcal{F}=\la X_1,\ldots,X_N\ra_{C^\infty_c(M)}$. Using Lemma \ref{lfgif},  $H_g\in N(\mathcal{I}_{\mathcal{F}})$ implies that for every $a=1,\ldots,N$ there exist functions $f_a^1,\ldots,f_a^N\in C^\infty(T^*M)$ such that

\begin{equation}\label{tay1}
     \tfrac{1}{2}\overline{\mathcal{L}_{X_a}g^{-1}}=\left \{
\overline{X_a},H_g\right \}=\sum_{b=1}^N f_a^b\overline{X_b},
\end{equation}
where we used Corollary \ref{symvect} for the first equality. Locally, for each $a$ and $b$, consider the first-order Taylor approximation of $f_a^b(q,p)$ with respect to the fiber coordinates $\{p_i\}$ around $(q,0)$

\begin{equation}\label{tay2}
f_a^b(q,p)=f_a^b(q,0)+\lambda_a^b(q,p)+o(\|p\|)
\end{equation}
where $\lambda_a^b(q,p)$ is linear in fiber coordinates. Equations \eqref{tay1} and \eqref{tay2} then imply that

\begin{equation*}
    \left \{
\overline{X_a},H_g\right \}-\sum_{b=1}^N\lambda_a^b\overline{X_b}=\sum_{b=1}^N\left[f_a^b(q,0)+o(\|p\|)\right]\overline{X_b}.
\end{equation*}
The left-hand side of the last equation is quadratic in fiber-coordinates, while the righ-hand side is not. This implies that both sides are identically zero. Consequently

\begin{align*}
    \left \{
\overline{X_a},H_g\right \}=\sum_{b=1}^N\lambda_a^b\overline{X_b},
\end{align*}
for some $\lambda_a^b\in C_1^\infty(T^*M)$. Lemma \ref{srfgen} then implies that $\mathcal{L}_Xg^{-1}\in \mathfrak{X}(M)\odot\mathcal{F}$. 
\end{proof}

\noindent Now we are able to present the proof of Proposition \ref{alggeo}.

\noindent\begin{proof}[{\textbf{Proposition} \ref{alggeo}}]
Let $(M,g,\mathcal{F})$ be a module SRF. As the statement is local, we can assume that $\mathcal{F}$ is finitely generated, i.e.\ there exist vector fields $X_1,...,X_N\in \mathfrak{X}(M)$ for some positive integer $N$, such that $\mathcal{F}=\left\langle X_1,...,X_N\right\rangle_{C^{\infty}_c(M)}$. By Lemma \ref{lfgif}, $\mathcal{I}_{\mathcal{F}}$ is generated by functions $\overline{X_1},...,\overline{X_N}$. By Proposition \ref{srfcons}, for every $a=1,\ldots,N$ there exist functions $\lambda_a^1,\ldots,\lambda_a^N\in C^\infty_1(T^*M)$ such that
\begin{equation*}
\left \{ H_g,\overline{X_a}\right \}=\sum_{b=1}^N\lambda_a^b\overline{X_b}.
\end{equation*}
Assume that $\gamma\colon(-\epsilon,\epsilon)\to M$ is a geodesic such that $\Dot{\gamma}(0)\perp F_{\gamma(0)}$, i.e.\ the geodesic is orthogonal to the leaf at $t=0$. Then the ideal $\mathcal{I}_{\mathcal{F}}(T^*M)$ vanishes at $(q_0,p_0)=(\gamma(0),g(\Dot{\gamma}(0),.))\in T^*M$. Since $\Phi_{H_g}^t(q_0)$ is defined for $-\epsilon<t<\epsilon$, for every $r<\epsilon$ there exists an open neighborhood $U\subset M$ of $q_0$ such that $\Phi_{H_g}^t$ is defined for $t\in(-r,r)$ on $U$.
According to  Proposition \ref{Hamflow},

\begin{align*}
g(\Dot{\gamma}(t), X_a(\gamma(t)))&=\overline{X_a}(\gamma(t),g_{\gamma(t)}(\Dot{\gamma}(t),.)) \\
&=\overline{X_a}\circ\Phi_{H_g}^t\vert_U(\gamma(0),g_{\gamma(0)}(\Dot{\gamma}(0),.))\\&=\Phi_{H_g}^t\vert_U^*\overline{X_a}(q_0,p_0) \, 
\end{align*} 
\noindent for $a=1,\ldots,N$. But now, according to Proposition \ref{normcons}, the function $(\Phi_{H_g}^t\vert_U)^*\left(\overline{X_a}\vert_{\Phi_{H_g}^t(U)}\right)$ is an element in $\mathcal{I}_{\mathcal{F}}(U)$ for $t\in(-r,r)$. This means that for $t$ in this interval, $\Dot{\gamma}(t)\perp F_{\gamma(t)}$. As $r<\epsilon$ is arbitrary, the proof is complete.
\end{proof}

\section{The functor \texorpdfstring{$\Phi$}{Phi} and reduction} \label{functor}

\noindent At the end of Section \ref{srfsection} we introduced the category $\mathbf{SRF}_0$ and in Section \ref{i-poisson} we introduced the category of $\mathcal{I}$-Poisson manifolds  \textbf{IPois}.
In this section we will provide a functor from the first to the second category, by sending a module SRF $(M,g,\mathcal{F})$ to the $\mathcal{I}$-Poisson manifold $(T^*M,\left\{\cdot,\cdot\right\},\mathcal{I}_{\mathcal{F}})$ and every surjective Riemannian submersion $\pi\colon(M_1,g_1)\to(M_2,g_2)$ to the bundle map $\varphi_{\pi}:=(g_2)_{\flat}\circ \exd\pi\circ (g_2)_{\flat}^{-1}\colon(T^*M_1,\left\{\cdot,\cdot\right\}_1)\to (T^*M_2,\left\{\cdot,\cdot\right\}_2)$, see Theorem \ref{thm:functor} below. $\varphi_{\pi}$ is precisely the map making the following diagram commute: 

\begin{center}
 \begin{tikzpicture}
  \node (A) {$TM_1$};
  \node (B) [below=of A] {$TM_2$};
  \node (C) [right=of A] {$T^*M_1$};
  \node (D) [right=of B] {$T^*M_2$};
  \draw[-stealth] (A)-- node[left] {$\exd\pi$} (B);
  \draw[-stealth] (B)-- node [below] {$(g_2)_{\flat}$} (D);
  \draw[-stealth] (A)-- node [above] {$(g_1)_{\flat}$} (C);
  \draw[-stealth] (C)-- node [right]{$\varphi_{\pi}$} (D);
\end{tikzpicture}
\end{center}
The map $\varphi_{\pi}$ is not a Poisson map in general:\footnote{In contrast to what is claimed in \cite{baez}.}

\begin{Example}
Let $\pi\colon\R^3\to \R^2$ be the projection $(x,y,z)\mapsto (x,y)$ in the canonical coordinates.  Equipping $\R^3$ and $\R^2$ with the following metric tensors

\begin{align*}
    g_3 :&= \exd x \otimes \exd x + (1+x^2)\exd y \otimes \exd y - x \exd y \otimes \exd z -x \exd z \otimes \exd y + \exd z \otimes \exd z \, \quad  \textit{and} \\
    g_2 :&= \exd \underline{x} \otimes \exd \underline{x} + \exd \underline{y} \otimes \exd \underline{y} , 
\end{align*}
respectively, turns $\pi$ into a Riemannian submersion.
Here $(\underline{x},\underline{y})$ denote the coordinates on $\R^2$. In the induced coordinates $(\underline{x},\underline{y},p_{\underline{x}},p_{\underline{y}})$ and $(x,y,z,p_x,p_y,p_z)$ on $T^*\R^2$ and $T^*\R^3$, respectively, the map $\varphi_{\pi}$ is given by
\[\varphi_{\pi}(x,y,z,p_x,p_y,p_z)=(x,y,p_x,p_y+xp_z).\]
This is not a Poisson map, since $\{\varphi_{\pi}^* p_{\underline{x}},\varphi_{\pi}^* p_{\underline{y}} \}=\{p_x,p_y+xp_z\}=p_z \neq 0$.
\end{Example}

\noindent In the last example the obstruction for $\varphi_{\pi}$ to be a Poisson map is that the horizontal distribution of the Riemannian submersion $\pi$, which is  generated by vector fields $\pderiv{}{x}$ and $\pderiv{}{y}+x\pderiv{}{z}$, is not integrable; the corresponding connection has  curvature.

\noindent The map $\varphi_{\pi}$ still preserves the Poisson bracket up to some ideal of functions $\mathcal{I}_{\ker \exd\pi}$:

\begin{Definition}
Let $\pi\colon M_1\to M_2$ be a submersion. The subsheaf of smooth functions  \emph{$\mathcal{I}_{\ker \exd\pi}$} on $T^*M_1$ is defined as the corresponding sheaf $\mathcal{I}_\mathcal{F}$ for the regular foliation $\mathcal{F}:=\Gamma_c(\ker \exd\pi)$.
\end{Definition}

\noindent When there is no ambiguity, for simplicity, we denote the ideal $\mathcal{I}_{\ker \exd\pi}(T^*M_1)$ by $\mathcal{I}_{\ker \exd\pi}$.

\begin{Remark}
It is not difficult to see that for every open subset $U\subset T^*M$, The ideal $\mathcal{I}_{\ker \exd\pi}(U)$ is the vanishing ideal of the submanifold $\mathrm{Ann}(\ker\exd\pi)\cap U\subset U$. Here $\mathrm{Ann}(\ker\exd\pi)$ stands for the annihilator of the subbundle $\ker\exd\pi\subset TM_1$. Moreover, since $\mathrm{Ann}(\ker\exd\pi)$ is an embedded submanifold, we have:

\begin{equation}\label{c}
    \mathcal{C}:=\left\{(q,p)\in T^*M_1\colon f(q,p)=0 \quad \forall f\in \mathcal{I}_{\ker \exd \pi}\right\}\equiv \mathrm{Ann}(\ker \exd \pi).
\end{equation}
\end{Remark}

\begin{Lemma}\label{sfpoisson}
Let $\pi\colon(M_1,g_1)\to(M_2,g_2)$ be a Riemannian submersion. Then for every $f,g\in C^{\infty}(T^*M_2)$:

\begin{align}
    \left\{f\circ \varphi_{\pi},g\circ \varphi_{\pi}\right\}_1-\left\{f,g\right\}_2\circ\varphi_{\pi}&\in& \mathcal{I}_{\ker \exd\pi}\label{i}\\
    \left\{f\circ\varphi_{\pi},\mathcal{I}_{\ker \exd\pi}\right\}_1&\subset& \mathcal{I}_{\ker \exd\pi}\label{ii}.
\end{align}
\end{Lemma}

\noindent \begin{proof} Choose local Darboux coordinates $(q_2^i,p_i^2)$ on $T^*M_2$ and $(q_1^i,q_1^{\alpha},p_i^1,p_{\alpha}^1)$ on $T^*M_1$, such that $q_2^i\circ \pi=q_1^i$ (this is possible since $\pi$ is assumed to be a submersion). In particular, $\mathcal{I}_{ker\exd \pi}$ is generated by the momenta $p_{\alpha}^1$. Now note that at every point $q\in M_1$, \[\exd_q\pi(\tfrac{\partial}{\partial q_1^i}|_{q})=\tfrac{\partial}{\partial q_2^i}|_{\pi(q)} , \]  since for every function $f\in C^{\infty}(M_2)$

\begin{equation*}
    \exd_q\pi(\tfrac{\partial}{\partial q_1^i}|_{q})\cdot f=\tfrac{\partial(f\circ \pi)}{\partial q_1^i}(q)=\tfrac{\partial f}{\partial q_2^j}(\pi(q))\tfrac{\partial(q_2^j\circ \pi)}{\partial q_1^i}(q)=\tfrac{\partial f}{\partial q_2^i}(\pi(q)).
\end{equation*}

\noindent In particular, since $\varphi_\pi$ is a bundle map, we have
\begin{equation} \label{qq}
   q_1^i =q_2^i\circ \varphi_{\pi}.
\end{equation}
Next we prove that upon restriction to  the vanishing submanifold $\mathcal{C}$ of $\mathcal{I}_{ker\exd \pi}$,

\begin{equation}
    \mathcal{C}=\left\{(q,p)\in T^*M_1\colon f(q,p)=0 \quad \forall f\in \mathcal{I}_{\ker \exd \pi}\right\}\equiv \mathrm{Ann}(\ker \exd \pi),
\end{equation}
one has $p_i^1=p_i^2\circ \varphi_{\pi}$. Indeed, let $(q,p)$ be a point in $T^*M_1$ and $X=(g_1)_{\flat}^{-1}(p)$. Then 

\begin{equation*}
    p_i^1(q,p)=p(\tfrac{\partial}{\partial q_1^i}|_q)=g_1(X,\tfrac{\partial}{\partial q_1^i}|_q)=g_1(X,(\tfrac{\partial}{\partial q_1^i}|_q)^H)+g_1(X,(\tfrac{\partial}{\partial q_1^i}|_q)^V)
\end{equation*}
\noindent where $(\tfrac{\partial}{\partial q_1^i}|_q)^H$ and $(\tfrac{\partial}{\partial q_1^i}|_q)^V$ are the horizontal and vertical parts of the vector $\tfrac{\partial}{\partial q_1^i}|_q$ with respect to $g_1$, respectively. Using that $\pi$ is a Riemannian submersion and that there exist functions $A_{\alpha}$ such that $(\tfrac{\partial}{\partial q_1^i}|_q)^V=\sum_\alpha A_{\alpha}(q)\tfrac{\partial}{\partial q_1^{\alpha}}|_q$, this implies:

\begin{equation*}
    p_i^1(q,p)=g_2(\exd_q\pi(X),\tfrac{\partial}{\partial q_2^i}|_{\pi(q)})+ \sum_\alpha g_1(X,A_{\alpha}(q)\tfrac{\partial}{\partial q_1^{\alpha}}|_q) .
\end{equation*}
 Consequently, by definition of $\varphi_{\pi}$,
 
\begin{equation*}
    p_i^1(q,p)=p_i^2\circ \varphi_{\pi}(q,p)+\sum_\alpha A_{\alpha}(q)p_{\alpha}^1(q,p)
\end{equation*}
 and thus 
 
\begin{equation} \label{pp}
   p_i^1\vert_\mathcal{C} =\left(p_i^2\circ \varphi_{\pi}\right)\!\vert_\mathcal{C}  .
\end{equation}
Now for every $f\in C^{\infty}(T^*M_2)$, we have:

\begin{equation}\label{chainq}
    \tfrac{\partial(f\circ \varphi_{\pi})}{\partial q_1^i}(q,p)=\tfrac{\partial f}{\partial q_2^j}(\varphi_\pi(q,p))\, \tfrac{\partial(q_2^j\circ \varphi_\pi)}{\partial q_1^i}(q,p)+\tfrac{\partial f}{\partial p_j^2}(\varphi_{\pi}(q,p))\,\tfrac{\partial(p_j^2\circ \varphi_{\pi})}{\partial q_1^i}(q,p).
\end{equation}

\noindent Since $\tfrac{\partial}{\partial q_1^i}\vert_\mathcal{C}$
is tangent to $\mathcal{C}$, for every point $(q,p)\in \mathcal{C}$, we may use Equation \eqref{pp}  to transform Equation \eqref{chainq} into:

\begin{align}\label{chainq2}
    \tfrac{\partial(f\circ \varphi_{\pi})}{\partial q_1^i}(q,p)&=\tfrac{\partial f}{\partial q_2^j}(\varphi_\pi(q,p))\tfrac{\partial q_1^j}{\partial q_1^i}(q,p)+\tfrac{\partial f}{\partial p_j^2}(\varphi_{\pi}(q,p))\tfrac{\partial p_j^1}{\partial q_1^i}(q,p)\\\notag
    &=\tfrac{\partial f}{\partial q_2^i}(\varphi_\pi(q,p)).
\end{align}

\noindent In a similar way, using the chain rule and that $\tfrac{\partial}{\partial p^1_i}\vert_\mathcal{C}$  and $\tfrac{\partial}{\partial q_1^{\alpha}}\vert_\mathcal{C}$  are tangent to $\mathcal{C}$, for every function $f\in C^{\infty}(T^*M_2)$ and every  $(q,p)\in \mathcal{C}$, one finds

\begin{align}\label{chainp}
\tfrac{\partial (f\circ \varphi_{\pi})}{\partial p_i^1}(q,p)&=\tfrac{\partial f}{\partial p_i^2}(\varphi_{\pi}(q,p)), \\
\tfrac{\partial(f\circ \varphi_{\pi})}{\partial q_1^\alpha}(q,p)&=0\label{chainalph}.
\end{align}

\noindent For every two functions $f,g\in C^{\infty}(T^*M_2)$, upon restriction to $\vert_\mathcal{C}$ we have:

\begin{align}\notag
    \left\{f\circ \varphi_{\pi}),g\circ \varphi_{\pi}\right\}_1\vert_\mathcal{C}&=\sum_{i}\left(\tfrac{\partial \left(f\circ \varphi_{\pi}\right)}{\partial p_i^1}\tfrac{\partial \left(g\circ \varphi_{\pi}\right)}{\partial q_1^i}-\tfrac{\partial \left(g\circ \varphi_{\pi}\right)}{\partial p_i^1}\tfrac{\partial \left(f\circ \varphi_{\pi}\right)}{\partial q_1^i}\right)\!{\vert_\mathcal{C}}\\\notag
    &+&\sum_{\alpha}\left(\tfrac{\partial \left(f\circ \varphi_{\pi}\right)}{\partial p_\alpha^1}\tfrac{\partial \left(g\circ \varphi_{\pi}\right)}{\partial q_1^{\alpha}}-\tfrac{\partial \left(g\circ \varphi_{\pi}\right)}{\partial p_\alpha^1}\tfrac{\partial \left(f\circ \varphi_{\pi}\right)}{\partial q_1^{\alpha}}\right)\!\vert_\mathcal{C}\\\notag
    &=\sum_{i}\left(\tfrac{\partial \left(f\circ \varphi_{\pi}\right)}{\partial p_i^1}\tfrac{\partial \left(g\circ \varphi_{\pi}\right)}{\partial q_1^i}-\tfrac{\partial \left(g\circ \varphi_{\pi}\right)}{\partial p_i^1}\tfrac{\partial \left(f\circ \varphi_{\pi}\right)}{\partial q_1^i}\right)\!{\vert_\mathcal{C}}\\\notag
    &=\sum_{i}\left(\tfrac{\partial f}{\partial p_i^2}\circ \varphi_{\pi}\right)\!{\vert_\mathcal{C}} \left( \tfrac{\partial g}{\partial q_2^i}\circ \varphi_{\pi}\right)\!{\vert_\mathcal{C}} -\left(\tfrac{\partial g}{\partial p_i^2}\circ \varphi_{\pi}\right)\!{\vert_\mathcal{C}} \left(\tfrac{\partial f}{\partial q_2^i}\circ \varphi_{\pi}\right)\!{\vert_\mathcal{C}}\\\notag
    &=\left\{f,g\right\}_2\circ \varphi_{\pi}{\vert_\mathcal{C}}
\end{align}
Here in the first equality we used just the definition of the Poisson bracket, in the second one we used Equation \eqref{chainalph}, 
thereafter Equations \eqref{chainq2} and \eqref{chainp}, and finally again the definition of the bracket.  Note that every function on $T^*M_1$ vanishing on $\mathcal{C}$ is an element of $\mathcal{I}_{\ker \exd\pi}$, which proves Equation \eqref{i}.

\noindent Equation \eqref{chainalph} implies Equation \eqref{ii} as well, since $\mathcal{I}_{\ker \exd\pi}$ is locally generated by coordinate functions $p_\alpha$ for $\alpha=1,\ldots,k$, and we have
\[\left\{p_\alpha,f\circ \varphi_{\pi}\right\}_1\vert_\mathcal{C}=\tfrac{\partial(f\circ \varphi_{\pi})}{\partial q_1^\alpha}\vert_\mathcal{C}=0.\] which gives $\left\{f\circ \varphi_{\pi},p_\alpha\right\}_1\in \mathcal{I}_{\ker \exd\pi}$.
\end{proof}

\begin{Corollary}\label{csubmn}
The restriction $\varphi_\pi\vert_\mathcal{C}\colon \mathcal{C}\to T^*M_2$ is a surjective submersion. It coincides with the projection to the leaf space for the coisotropic reduction of $\mathcal{C}\subset T^*M_1$.
\end{Corollary}

\noindent \begin{proof} Choosing the same local coordinates as in the proof of Lemma \ref{sfpoisson}, $\left(q_1^i,q_1^{\alpha},p_i^1\right)$ give local coordinates for $\mathcal{C}$ and Equations \eqref{qq} and \eqref{pp} ensure that $\varphi_\pi\left(q_1^i,q_1^{\alpha},p_i^1\right)=\left(q_1^i,p_i^1\right)$.
\end{proof}

\noindent To study the obstruction for $\varphi_{\pi}$ to be a Poisson map, we first prove the following lemma which describes the horizontal distribution in terms of the map $\varphi_{\pi}$:

\begin{Lemma}\label{horpoiss}
Let $\pi \colon (N,h)\to (M,g)$ be a Riemannian submersion and let $X$ be a vector field on $M$. Then the horizontal lift of $X$ is given by a vector field $V$ on $N$ satisfying

\begin{equation}
\overline{V}=(\varphi_{\pi})^*\overline{X} ,
\end{equation}
which is an element in $C_1^{\infty}(T^*N)$.
\end{Lemma}

\noindent\begin{proof}
Define $V\in\mathfrak{X}(N)$ by $\overline{V}:=(\varphi_{\pi})^*(\overline{X})\in C_1^{\infty}(T^*M)$. Using Corollary \ref{symvect} and Lemma \ref{sfpoisson} for every function $f\in C^{\infty}(M)$, we have

\begin{align*}
    V\cdot \pi^*f&=\left\{\overline{V},\overline{\pi^*f}\right\}=\left\{(\varphi_{\pi})^*\overline{X},\pi^*f\right\}\\
                   &=\left\{(\varphi_{\pi})^*\overline{X},(\varphi_{\pi})^*f\right\}=(\varphi_{\pi})^*\left\{\overline{X},f\right\}\\
                   &=\pi^*(X\cdot f),
\end{align*}
which means that $V$ is projectable and projects to $X$. In addition, for every vertical vector $v\in \ker \exd_q\pi$ and $i=1,\ldots,n$ we have
\begin{align*}
    h(v,V\vert_q)&=\overline{V}(h_{\flat}(v))=(\varphi_{\pi})^*\overline{X}(h_{\flat}(v))\\
            &=\overline{X}(\varphi_{\pi}\circ h_{\flat}(v))=\overline{X}(g_{\flat}\circ\exd_q \pi (v))\\
            &=0,
\end{align*}
\noindent showing that $V$ is the horizontal lift of $X$.
\end{proof} 

\noindent The following identifies the obstruction for $\varphi_{\pi}$ to be a Poisson map:

\begin{Proposition}
Let $\pi\colon(N,h)\to (M,g)$ be a Riemannian submersion. Then the map $\varphi_{\pi}:=g_{\flat}\circ \exd\pi\circ h_{\flat}^{-1}$ is a Poisson map if and only if the horizontal  distribution $\mathcal{H}\subset TM$ of $\pi$ is integrable.
\end{Proposition}

\noindent \begin{proof}
Let $\left\{f_i\right\}_{i=1}^n$ be a local orthonormal frame around a point $\underline{q}\in M$ and $\left\{e_i\right\}_{i=1}^n$ their horizontal lifts. By Lemma \ref{horpoiss} we have $\overline{e_i}=(\varphi_{\pi})^*\overline{f_i}$ for $i=1,\ldots,n$. If $\varphi_{\pi}$ is a Poisson map, the family of functions $(\varphi_{\pi})^*(\overline{f_i})\in C_1^{\infty}(T^*N)$ is closed under the Poisson bracket, and consequently the horizontal distribution locally generated by vector fields $e_i$ is integrable. This proves the if part of the proposition.

\noindent Conversely assume that $\mathcal{H}$ is integrable. It is enough to check the condition of being a Poisson map on smooth functions in $C^\infty_0(T^*M)\bigoplus C^\infty_1(T^*M)$ only. First, for every $f,g\in C^\infty_0(T^*M)$ we have $\left\{f\circ\varphi_\pi,g\circ\varphi_\pi\right\}_{T^*N}=\left\{f,g\right\}_{T^*M}=0$. Second,  for every $\overline{X}\in C^\infty_1(T^*M)$ and $f\in C^\infty_0(T^*M)$ we have

\begin{align*}
    \left\{\overline{X}\circ\varphi_\pi,f\circ\varphi_\pi\right\}_{T^*N}&=X^H\cdot(f\circ\varphi_\pi)\\
    &=(X\cdot f)\circ \varphi_\pi\\&=\left\{\overline{X},f\right\}_{T^*M}\circ\varphi_\pi.
\end{align*}
Finally, by Lemma \ref{horpoiss} and integrability of $\mathcal{H}$, for every $\overline{X},\overline{Y}\in C^\infty_1(T^*M)$ one obtains
\begin{align*}
    \left\{\overline{X}\circ\varphi_\pi,\overline{Y}\circ\varphi_\pi\right\}_{T^*N}&=\overline{[X^H,Y^H]}\\
    &=\overline{[X,Y]^H}\\&=\overline{[X,Y]}\circ \varphi_\pi\\&=\left\{\overline{X},\overline{Y}\right\}_{T^*M}\circ\varphi_\pi.
\end{align*}
\end{proof} 

\begin{Lemma}\label{metricdiff}
Let $\pi\colon(N,h)\to (M,g)$ be a Riemannian submersion. Then
\[H_h-H_g\circ \varphi_{\pi}\in \mathcal{I}_{\ker \exd\pi}.\]
\end{Lemma}

\noindent\begin{proof}
It is enough to show that the left-hand side vanishes on $\mathcal{C}$, defined in Equation \eqref{c}. For every $(q,p)\in \mathcal{C}$, we have
\begin{align*}
    H_h(q,p)&=\tfrac{1}{2}\left\la p,(h_{\flat})^{-1}(p)\right\ra=\tfrac{1}{2}\left\la p,[(h_{\flat})^{-1}(p)]^H\right\ra\\
    &=\tfrac{1}{2}h((h_{\flat})^{-1}(p),[(h_{\flat})^{-1}(p)]^H)\\
    &=\tfrac{1}{2}h([(h_{\flat})^{-1}(p)]^H,[(h_{\flat})^{-1}(p)]^H)\\
    &=\tfrac{1}{2}g(\exd_q\pi\circ(h_{\flat})^{-1}(p),\exd_q\pi\circ(h_{\flat})^{-1}(p))\\
    &=\tfrac{1}{2}\left\la \varphi_{\pi}(p),(g_{\flat})^{-1}(\varphi_{\pi}(p))\right\ra\\
    &=H_g\circ \varphi_{\pi}(q,p).
\end{align*}
\end{proof} 

\noindent Now we are able to prove the well-behavedness of module SRFs under Riemannian submersions.

\noindent \begin{proof}[{\textbf{Proposition} \ref{pullback}}]
Let $(M,g,\mathcal{F})$ be a module SRF and $\pi\colon(N,h)\to (M,g)$ a Riemannian submersion. By Lemma \ref{decomposition}

\begin{equation}\label{split}
    \mathcal{F}_N=\langle\mathcal{F}_M^{\mathcal{H}}+\Gamma\left(\ker \exd\pi\right) \rangle_{C^{\infty}_c(N)}
\end{equation}
\noindent where $\mathcal{F}_M^{\mathcal{H}}$ consists of horizontal lifts of vector fields in $\mathcal{F}_N$. By Lemma \ref{horpoiss} we have

\begin{equation*}
    \mathcal{J}_{\mathcal{F}_N}=\left\langle(\varphi_{\pi})^*\overline{\mathcal{F}_M}+\mathcal{I}_{\ker \exd\pi}\right\rangle_{C^{\infty}_c(T^*N)},
\end{equation*}
where $\overline{\mathcal{F}_M}:=\left\{\overline{X}\,\colon\,X\in \mathcal{F}_M\right\}$. Finally, it remains to check $H_g\in N(\mathcal{I}_{\pi^{-1}\mathcal{F}})$. By Lemma \ref{normif} it is enough to verify the following: 

\begin{align*}
    \left\{H_g,(\varphi_{\pi})^*\overline{\mathcal{F}_M}+\mathcal{I}_{\ker \exd\pi}\right\}&=\left\{H_g-H_h\circ \varphi_{\pi},(\varphi_{\pi})^*\overline{\mathcal{F}_M}\right\}+\left\{H_h\circ \varphi_{\pi},(\varphi_{\pi})^*\overline{\mathcal{F}_M}\right\}\\
    &+&\left\{H_g-H_h\circ \varphi_{\pi},\mathcal{I}_{\ker \exd\pi}\right\}+\left\{H_h\circ \varphi_{\pi},\mathcal{I}_{\ker \exd\pi}\right\}\\
    &\subset&(\varphi_{\pi})^*\overline{\mathcal{F}_M}+\mathcal{I}_{\ker \exd\pi}.
\end{align*}
Here we used Lemmas \ref{sfpoisson} and \ref{metricdiff} to prove the inclusion.
\end{proof}

\noindent The following  theorem is the main result of this section.

\begin{theorem}\label{thm:functor}
The map sending every module SRF $(M,g,\mathcal{F})$ to the corresponding dynamical $\mathcal{I}$-Poisson manifold $(T^*M,\left\{\cdot,\cdot\right\}_{T^*M},\mathcal{I}_\mathcal{F},H_g)$ and every morphism $\pi$ of SRFs within $\mathbf{SRF}_0$ to the map $\varphi_{\pi}$ defines a functor  $\Phi\colon\mathbf{SRF}_0\to \mathbf{dynIPois}$.
\end{theorem}

\noindent \begin{proof} It is enough to show that $\Phi$ preserves the morphisms. A morphism $\pi$ within $\mathbf{SRF}_0$ is a Riemannian submersion $\pi\colon(N,h,\mathcal{F}_N)\to (M,g,\mathcal{F}_M)$ such that $\mathcal{F}_N=\pi^{-1}\left(\mathcal{F}_M\right)$. Similar to the previous proof we have

\begin{equation*}
   \mathcal{J}_{\mathcal{F}_N}=\left\langle(\varphi_{\pi})^*\overline{\mathcal{F}_M}+\mathcal{I}_{\ker \exd\pi}\right\rangle_{C^{\infty}_c(T^*N)}\,
\end{equation*}
and therefore the pullback $(\varphi_{\pi})^*\mathcal{I}_{\mathcal{F}_M}(T^*M)$ and $\mathcal{I}_{\ker \exd\pi}$ lie inside $\mathcal{I}_{\mathcal{F}_N}(T^*N)$. By Lemma \ref{sfpoisson}, for every $f\in N(\mathcal{I}_{\mathcal{F}_M})$ we have

\begin{align*}
    \left\{f\circ \varphi_{\pi},(\varphi_{\pi})^*\overline{\mathcal{F}_M}+\mathcal{I}_{\ker \exd\pi}\right\}_{T^*N}&\subset&\left\{f,\overline{\mathcal{F}_M}\right\}_{T^*M}\circ \varphi_{\pi}+\left\{f\circ \varphi_{\pi},\mathcal{I}_{\ker \exd\pi}\right\}_{T^*N}+\mathcal{I}_{\ker \exd\pi}\\
    &\subset& (\varphi_{\pi})^*\overline{\mathcal{F}_M}+\mathcal{I}_{\ker \exd\pi}
\end{align*}
which implies that $\left\{f\circ \varphi_{\pi},\mathcal{J}_{\mathcal{F}_N}\right\}_{T^*N}\subset \mathcal{J}_{\mathcal{F}_N}$, and consequently $(\varphi_{\pi})^*N(\mathcal{I}_{\mathcal{F}_M})$ lies inside $N(\mathcal{I}_{\mathcal{F}_N})$. Using Lemma \ref{sfpoisson} again, for every $f,g\in N(\mathcal{I}_{\mathcal{F}_M})$ 

\begin{align*}
     \left\{f\circ \varphi_{\pi},g\circ \varphi_{\pi}\right\}_{T^*N}-\left\{f,g\right\}_{T^*M}\circ \varphi_{\pi}\in \mathcal{I}_{\ker \exd\pi}\subset \mathcal{I}_{\mathcal{F}_N}(T^*N).
\end{align*}
These together with Lemma \ref{metricdiff} complete the proof.
\end{proof}

 \begin{theorem}\label{thmc}
Let $(M_1,\mathcal{F}_1)$ and $(M_1,\mathcal{F}_1)$ be Hausdorff Morita equivalent singular foliations. Then the Poisson algebras $\mathcal{R}\left({\mathcal{I}_{\mathcal{F}_1}}\right)$ and  $\mathcal{R}\left({\mathcal{I}_{\mathcal{F}_2}}\right)$
are isomorphic.
\end{theorem}

\noindent Here $\mathcal{R}\left({\mathcal{I}_{\mathcal{F}_i}}\right) \equiv N(\mathcal{I}_{\mathcal{F}_i})/\mathcal{I}_{\mathcal{F}_i}(T^*M_i)$, $i=1,2$, see Definition \ref{reducedPois}.

\noindent The proof of this theorem will be a consequence of the following two lemmas. 

\begin{Lemma}\label{inj}
Let $\pi\colon (N,h)\to (M,g)$ be a surjective Riemannian submersion with connected fibers and $\mathcal{F}$ be an SF on $M$. If $f\circ\varphi_\pi\in\mathcal{I}_{\pi^{-1}\mathcal{F}}(T^*N)$ for some $f\in C^\infty(T^*M)$, then $f\in \mathcal{I}_{\mathcal{F}}(T^*M)$.
\end{Lemma}

\noindent \begin{proof}
We first demonstrate that the result holds true for finitely generated singular foliations. Let $\mathcal{F}=\la X_1,\ldots,X_N\ra_{C^\infty_c(M)}$ and let $Y_1,\ldots,Y_K\in \mathfrak{X}(N)$ be generators of the regular foliation $\Gamma_c(\ker\exd\pi)$ for some positive integers $N$ and $K$. Lemma \ref{decomposition} implies that $\pi^{-1}\mathcal{F}=\langle X_1^{\mathcal{H}},\ldots,X_N^{\mathcal{H}},Y_1,\ldots,Y_K \rangle_{C^{\infty}_c(N)}$. Consequently, for every open subset $V\subset T^*N$, we obtain 

\begin{equation}\label{ip-f}
    \mathcal{I}_{\pi^{-1}\mathcal{F}}(V)=\la\overline{X_1}\circ \varphi_\pi\vert_V,\ldots,\overline{X_N}\circ \varphi_\pi\vert_V,\overline{Y_1}\vert_V,\ldots,\overline{Y_K}\vert_V\ra_{C^\infty(V)},
\end{equation}
where we used Lemmas \ref{horpoiss} and \ref{lfgif}.

\noindent Let us assume for a moment that there exists a global section $s\colon T^*M\to \mathcal{C}$ for the surjective submersion $\varphi_\pi\vert_\mathcal{C}\colon \mathcal{C}\to T^*M$ (see Corollary \ref{csubmn}). Since $f\circ\varphi_\pi\in\mathcal{I}_{\pi^{-1}\mathcal{F}}(T^*N)$, Equation \eqref{ip-f} implies that there exist smooth functions $\lambda^1,\ldots,\lambda^N,\eta^1,\ldots,\eta^K\in C^\infty(T^*N)$ such that

\begin{equation}\label{fophi}
    f\circ\varphi_\pi=\sum_{a=1}^N\lambda^a\cdot\left(\overline{X_a}\circ \varphi_\pi\right)+\sum_{b=1}^K\eta^b\cdot\overline{Y_b}.
\end{equation}
Since $\varphi_\pi\circ s\circ\varphi_\pi=\varphi_\pi$ and $\overline{Y_b}\vert_\mathcal{C}=0$, composing both sides of Equation \eqref{fophi} by $s\circ\varphi_\pi$ gives

\begin{align*}
    f\circ\varphi_\pi=f\circ\varphi_\pi\circ s\circ\varphi_\pi&=\sum_{a=1}^N\left(\lambda^a\circ s\circ\varphi_\pi\right)\cdot\left(\overline{X_a}\circ \varphi_\pi\circ s\circ\varphi_\pi\right)\\
    &=\left(\sum_{a=1}^N\left(\lambda^a\circ s\right)\cdot\overline{X_a}\right)\circ \varphi_\pi.
\end{align*}
This implies that $f=\sum_{a=1}^N\left(\lambda^a\circ s\right)\cdot\overline{X_a}\in \mathcal{I}_\mathcal{F}(T^*M)$, since $\varphi_\pi$ is surjective.

\noindent If a global section does not exist, we can choose an open covering $\left\{U_i\right \}_{i=1}^\infty$ of $T^*M$ such that for every positive integer $i$ there exists a local section $s_i\colon U_i\to \mathcal{C}$. Using the same argument as for the global case, we may show that $f\vert_{U_i}\in \mathcal{I}_\mathcal{F}(U_i)$ for each $i$. Since $\mathcal{I}_\mathcal{F}$ is a sheaf on $T^*M$, we have $f\in \mathcal{I}_\mathcal{F}(T^*M)$.

\noindent For the general case, choose an open covering $\left\{U_i\right \}_{i=1}^\infty$ of $M$ such that for every positive integer $i$ the pullback  $\iota_{U_i}^{-1}\mathcal{F}$ is finitely generated. The finitely generated case discussed before then implies that $f\vert_{T^*U_i}\in \mathcal{I}_{\mathcal{F}}(T^*U_i)$ for every $i$, and since $\mathcal{I}_{\mathcal{F}}$ is a sheaf, we obtain $f\in \mathcal{I}_{\mathcal{F}}(T^*M)$.
\end{proof}

\begin{Lemma}\label{surj}
Let $\pi\colon (N,h)\to (M,g)$ be a surjective Riemannian submersion with connected fibers and $\mathcal{F}$ be a finitely generated SF on $M$. Then for every $F\in N(\mathcal{I}_{\pi^{-1}\mathcal{F}})$, there exists some $f\in C^\infty(T^*M)$ such that $F-f\circ \varphi_\pi \in \mathcal{I}_{\pi^{-1}\mathcal{F}}(T^*N)$.
\end{Lemma}

\noindent \begin{proof}
We proceed as in the beginning of the proof of Lemma \ref{inj}, establishing Equation \eqref{ip-f} and assuming  first again that there is a global section $s\colon T^*M\to \mathcal{C}= \mathrm{Ann}(\ker \exd \pi)$ for the surjection $\varphi_\pi\vert_\mathcal{C}\colon \mathcal{C}\to T^*M$. Define, in addition, $f:=F\circ s \in C^\infty(T^*M)$. 

\noindent We now will prove that for every $x\in T^*M$, there exists an open neighborhood $V_x\subset T^*N$ such that $(F-f\circ \varphi_\pi)\vert_{V_x}\in  \mathcal{I}_{\pi^{-1}\mathcal{F}}(V_x)$. Since $\mathcal{I}_{\pi^{-1}\mathcal{F}}$ is a sheaf, this implies the desired $F-f\circ \varphi_\pi \in \mathcal{I}_{\pi^{-1}\mathcal{F}}(T^*N)$. The proof is divided into the following three cases:

\noindent \textbf{Case 1. [$x\not\in\mathcal{C}$]:} We choose an open subset $V_x\subset T^*N$ such that $\overline{V_x}\cap \mathcal{C}=\phi$. Let $\rho\in C^\infty(T^*N)$ with $\mathrm{supp}(\rho)\subset T^*N\setminus\mathcal{C}$ and $\rho\vert_{\overline{V_x}}\equiv 1$. Then  since $\rho\cdot (F-f\circ \varphi_\pi)$ vanishes on $\mathcal{C}$, we obtain
\[(F-f\circ \varphi_\pi)\vert_{V_x}=\rho\cdot (F-f\circ \varphi_\pi)\vert_{V_x}\in \mathcal{I}_{\ker \exd\pi}(V_x)\subset \mathcal{I}_{\pi^{-1}\mathcal{F}}(V_x).\]

\noindent \textbf{Case 2. [$x\in s(T^*M)\subset \mathcal{C}$]:} Choose local coordinates $(q^i,q^\alpha)$ centered at the base-point of $x$ and $(q^i_M)$ on $N$ and $M$, respectively, were $i\in \{1,\ldots,m := \dim M\}$ and $\alpha \in \{m+1,\ldots,n:= \dim N\}$, which are 
compatible with the submersion $\pi$, i.e.\ $\pi(q^i,q^\alpha)=(q^i)$.  Let $(q^i,q^\alpha,p_i,p_\alpha)$ be the corresponding Darboux coordinates on some open neighborhood $V_x\subset T^*N$ centered at $x$. As a consequence, in particular  $\varphi_\pi (q^i,q^\alpha,p_i,0) = (q^i,p_i)$ (see Corollary \ref{csubmn}) and $\mathcal{I}_{\ker \exd\pi}(V_x)=\la p_\alpha\ra_{C^\infty(V_x)}$. For simplicity also assume that, in these local coordinates, $s\circ \varphi_\pi (q^i,q^\alpha,p_i,0)=(q^i,0,p_i,0)$. Then, for every arbitrary point $(q_0^i,q_0^\alpha,p^0_i,0)$ in $V_x\cap \mathcal{C}$, we have
\begin{align}
    (F-f\circ \varphi_\pi)(q_0^i,q_0^\alpha,p^0_i,0)&=F(q_0^i,q_0^\alpha,p^0_i,0)-F(q_0^i,0,p^0_i,0)\notag\\
    &=\int_0^1\tfrac{\exd}{\exd t}F(q_0^i,tq_0^\alpha,p^0_i,0)\,\exd t\notag\\
    &=\int_0^1\left(\sum_\beta \left\{q_0^\beta p_\beta,F\right\}_{T^*N}\right)\!(q_0^i,tq_0^\alpha,p^0_i,0)\:\exd t\notag\\
    &=\int_0^1\left(\sum_\beta q_0^\beta \left\{p_\beta,F\right\}_{T^*N}\right)\!(q_0^i,tq_0^\alpha,p^0_i,0)\:\exd t .\label{case2}
\end{align}
Since $F\in N(\mathcal{I}_{\pi^{-1}\mathcal{F}})$, for every $\beta$ there exist smooth functions $\lambda_\beta^1,\ldots,\lambda_\beta^N,\eta_\beta^1,\ldots,\eta_\beta^K\in C^\infty(V_x)$ such that $\left\{p_\beta,F\right\}_{T^*N}=\sum_a\lambda_\beta^a\cdot (\overline{X_a}\circ \varphi_\pi\vert_{V_x})+\sum_b\eta_\beta^b\cdot(\overline{Y_b}\vert_{V_x})$. Implementing this into Equation \eqref{case2} gives

\begin{align}
    &\quad(F-f\circ \varphi_\pi)(q_0^i,q_0^\alpha,p^0_i,0)\notag\\
    &=\int_0^1\left(\sum_\beta q_0^\beta\left(\sum_a\lambda_\beta^a\cdot (\overline{X_a}\circ \varphi_\pi\vert_{V_x})\right)\right)\!(q_0^i,tq_0^\alpha,p^0_i,0)\:\exd t\notag\\
    &=\sum_a\left(\Lambda^a\cdot (\overline{X_a}\circ \varphi_\pi\vert_{V_x})\right)(q_0^i,q_0^\alpha,p^0_i,0)\label{case22}\, 
\end{align}
where $\Lambda^a\in C^\infty(V_x)$ is defined as
\[\Lambda^a(q^i,q^\alpha,p_i,p_\alpha):=\int_0^1\left(\sum_\beta q^\beta\lambda_\beta^a\right)\!(q^i,tq^\alpha,p_i,p_\alpha)\:\exd t.\]
Equation \eqref{case22} implies that $(F-f\circ\varphi_\pi)-\sum_a\Lambda^a\cdot (\overline{X_a}\circ \varphi_\pi\vert_{V_x})$ vanishes on $V_x\cap \mathcal{C}$ and consequently this difference is an element of $\mathcal{I}_{\ker \exd\pi}(V_x)$. Since $\sum_a\Lambda^a\cdot (\overline{X_a}\circ \varphi_\pi\vert_{V_x})\in \mathcal{I}_{\pi^{-1}\mathcal{F}}(V_x)$ we obtain
\[(F-f\circ\varphi_\pi)\vert_{V_x}\in \mathcal{I}_{\pi^{-1}\mathcal{F}}(V_x).\]

\noindent \textbf{Case 3. [$x\in \mathcal{C} \setminus s(T^*M)$]:} Define $x_0=s\circ\varphi_\pi(x)$. Since $x$ and $x_0$ belong to the same fiber of $\varphi_\pi\vert_\mathcal{C}$, there exist compactly supported functions $h_1,\ldots,h_l\in \mathcal{I}_{\ker \exd\pi}$ for some positive integer $l$ such that their Hamiltonian flows connect $x_0$ to $x$, i.e.\

\begin{align*}
    x=\Phi_{h_l}^1\circ\ldots\circ\Phi_{h_1}^1(x_0).
\end{align*}
Then the global section $s':=\Phi_{h_l}^1\circ\ldots\circ\Phi_{h_1}^1\circ s$ passes through the point $x$. After \textbf{Case 2} for the function $f':=F\circ s'$, there exists an open neighborhood $V_x$ around $x$ such that $(F-f'\circ\varphi_\pi)\vert_{V_x}\in \mathcal{I}_{\pi^{-1}\mathcal{F}}(V_x)$. It remains to show that $(f'-f)\circ \varphi_\pi\in \mathcal{I}_{\pi^{-1}\mathcal{F}}(T^*N)$. For arbitrary $y\in \mathcal{C}$, defining $\Phi_{h_0}^t:=\mathrm{Id}_{T^*N}$ and $y_0=s\circ \varphi_\pi(y)$  gives

\begin{align}
    (f'-f)\circ \varphi_\pi(y)&=(F\circ s'-F\circ s)\circ\varphi_\pi(y)\notag\\
    &=\sum_{i=1}^l(F\circ\Phi_{h_i}^1\circ\ldots\circ\Phi_{h_0}^1-F\circ\Phi_{h_{i-1}}^1\circ\ldots\circ\Phi_{h_0}^1)(y_0)\notag\\
    &=\sum_{i=1}^l\int_0^1 \tfrac{\exd}{\exd t}F\circ\Phi_{h_i}^t\circ\Phi_{h_{i-1}}^1\circ\ldots\circ\Phi_{h_0}^1(y_0)\:\exd t\notag\\
    &=\sum_{i=1}^l\int_0^1 \left\{h_i,F\right\}_{T^*N}\circ\Phi_{h_i}^t\circ\Phi_{h_{i-1}}^1\circ\ldots\circ\Phi_{h_0}^1(y_0)\:\exd t\label{case3}.
\end{align}
Since $F\in N(\mathcal{I}_{\pi^{-1}\mathcal{F}})$, for every $i$ there exist smooth functions $\lambda_i^1,\ldots,\lambda_i^N,\eta_i^1,\ldots,\eta_i^K\in C^\infty(T^*N)$ such that $\left\{h_i,F\right\}_{T^*N}=\sum_a\lambda_i^a\cdot (\overline{X_a}\circ \varphi_\pi)+\sum_b\eta_i^b\cdot(\overline{Y_b})$. Implementing this into Equation \eqref{case3}, making use of the fact that the flows of the $h_i$s preserve $\mathcal{C}$, and noting that the $\overline {Y_b}$s vanish on $\mathcal{C}$, this gives

\begin{align}
    &\quad(f'-f)\circ \varphi_\pi(y)\notag\\
    &=\sum_{i=1}^l\int_0^1 \left(\sum_a\lambda_i^a\cdot (\overline{X_a}\circ \varphi_\pi)\right)\circ\Phi_{h_i}^t\circ\Phi_{h_{i-1}}^1\circ\ldots\circ\Phi_{h_1}^1(y_0)\:\exd t\notag\\
    &=\sum_a\Lambda^a(y)\cdot (\overline{X_a}\circ \varphi_\pi(y))\label{case33}.
\end{align}
Here we defined $\Lambda^a\in C^\infty(T^*N)$ by 
\[\Lambda^a(z):=\sum_{i=1}^l\int_0^1 \lambda_i^a\circ\Phi_{h_i}^t\circ\Phi_{h_{i-1}}^1\circ\ldots\circ\Phi_{h_1}^1\circ s\circ\varphi_\pi(z)\:\exd t\qquad \forall\,z\in T^*N.\]
Equation \eqref{case33} implies that $ (f'-f)\circ \varphi_\pi-\sum_a\Lambda^a(y)\cdot (\overline{X_a}\circ \varphi_\pi(y))$ vanishes on $\mathcal{C}$ and, equivalently, it thus belongs to $\mathcal{I}_{\ker \exd\pi}$ and since $\sum_a\Lambda^a\cdot (\overline{X_a}\circ \varphi_\pi)\in \mathcal{I}_{\pi^{-1}\mathcal{F}}(T^*N)$. This gives  $(f'-f)\circ \varphi_\pi\in \mathcal{I}_{\pi^{-1}\mathcal{F}}(T^*N)$, which completes the proof in \textbf{Case 3}.

\noindent If a global section does not exist, we can choose a locally finite open covering $\left\{U_i\right \}_{i=1}^\infty$ of $T^*M$ with a partition of unity $\left\{\rho_i\right \}_{i=1}^\infty$ subordinate to it, such that for every positive integer $i$ there exists a local section $s_i\colon U_i\to \mathcal{C}$. Similar to the global case, we can show that for $f_i:=F\circ s_i\in C^\infty(U_i)$, we have $F\vert_{\varphi_\pi^{-1}(U_i)}-f\circ \varphi_\pi\vert_{\varphi_\pi^{-1}(U_i)}\in \mathcal{I}_{\pi^{-1}\mathcal{F}}(\varphi_\pi^{-1}(U_i))$. Defining $f:=\sum_{i=1}^\infty\rho_i f_i$, we claim that $F-f\circ\varphi_\pi\in \mathcal{I}_{\pi^{-1}\mathcal{F}}(T^*N)$. This is equivalent to showing that for every $\sigma\in C^\infty_c(T^*N)$ we have $\sigma\cdot(F-f\circ\varphi_\pi)\in \mathcal{J}_{\pi^{-1}\mathcal{F}}$. Since $\mathrm{supp}(\sigma)$ is compact, it can be covered by finitely many open subsets $\varphi_\pi^{-1}(U_{i_1}),\ldots,\varphi_\pi^{-1}(U_{i_n})$ of the covering $\left\{\varphi_\pi^{-1}(U_i)\right \}_{i=1}^\infty$. This gives

\begin{align}
    \sigma\cdot(F-f\circ\varphi_\pi)&=\sum_{i=1}^\infty\sigma\cdot(\rho_i\circ\varphi_\pi)\cdot(F-f\circ\varphi_\pi)\label{surjlast}\\&=\sum_{a=1}^n\sigma\cdot(\rho_{i_a}\circ\varphi_\pi)\cdot(F\vert_{\varphi_\pi^{-1}(U_{i_a})}-f_{i_a}\circ\varphi_\pi\vert_{\varphi_\pi^{-1}(U_{i_a})})\notag\\&\in& \mathcal{J}_{\pi^{-1}\mathcal{F}}\notag,
\end{align}
since $\sigma\cdot(\rho_{i_a}\circ\varphi_\pi)\in C^\infty_c(\varphi_\pi^{-1}(U_{i_a}))$. This completes the proof.
\end{proof}

\noindent \begin{proof}[{\textbf{Theorem} \ref{thmc}}]
It is enough to show that for every surjective submersion $\pi \colon N\to M$ with connected  fibers over a foliated manifold $(M,\mathcal{F})$, the Poisson algebras $\mathcal{R}\left({\mathcal{I}_{\mathcal{F}}}\right)$ and $\mathcal{R}\left(\mathcal{I}_{\pi^{-1}\mathcal{F}}\right)$ are isomorphic. 

\noindent To do so, we first choose Riemannian metrics $g_M$ and $g_N$ such that $\pi$ becomes a Riemannian submersion. This can be done as follows: choose a Riemannian metric $g_M$ on $M$, a fiber metric $g^\perp$ on $\ker \exd\pi\subset TN$, and a subbundle $\mathcal{H}\subset TN$ complementary to $\ker \exd\pi$; one then declares these two subbundles to be  orthogonal to one another and  defines  $g_N=\left(\pi^*g_M\right)\!\vert_\mathcal{H}+g^\perp$.

\noindent Injectivity of $\tilde{\varphi_\pi}$ is a direct consequence of Lemma \eqref{inj}. It remains to prove that $\tilde{\varphi_\pi}$ is surjective. It follows from showing that, for every $F\in N(\mathcal{I}_{\pi^{-1}\mathcal{F}})$,  there exists $f\in N(\mathcal{I}_{\mathcal{F}})$ such that $F-f\circ \varphi_\pi\in \mathcal{I}_{\pi^{-1}\mathcal{F}}(T^*N)$. To do so we 
choose an open covering $\left\{U_i\right\}_{i=1}^\infty$ of $M$ such that, for every positive integer $i$, the pullback  $\iota_{U_i}^{-1}\mathcal{F}$ is finitely generated. Let $\left\{V_a\right\}_{a=1}^\infty$ be a locally finite refinement  of the covering $\left\{T^*U_i)\right\}_{i=1}^\infty$ of $T^*M$ and let $\left\{\rho_a\right\}_{a=1}^\infty$ be a partition of unity  subordinate to $\left\{V_a\right\}_{a=1}^\infty$. Lemma \ref{surj} then implies that for every $a$ there exists $f_a\in C^\infty(V_a)$ such that $F\vert_{\varphi_\pi^{-1}(V_a)}-f_a\circ\varphi_\pi\vert_{\varphi_\pi^{-1}(V_a)}\in \mathcal{I}_{\pi^{-1}\mathcal{F}}(\varphi_\pi^{-1}(V_a))$. Using the same argument as in the proof of Lemma \ref{surj} (see Equation \eqref{surjlast}), for $f:=\sum_{a=1}^\infty\rho_af_a$ we have $F-f\circ\varphi_\pi\in \mathcal{I}_{\pi^{-1}\mathcal{F}}(T^*N)$. To complete the proof, we show that $f\in N(\mathcal{I}_\mathcal{F})$ as follows: Since $f\circ\varphi_\pi\in N(\mathcal{I}_{\pi^{-1}\mathcal{F}})$, Equation \eqref{i} of Lemma \ref{sfpoisson} implies that $\left\{f,\mathcal{J}_{\mathcal{F}}\right\}_{T^*M}\circ\varphi_\pi\subset\varphi_\pi^*\mathcal{J}_{\mathcal{F}}\subset \mathcal{I}_{\pi^{-1}\mathcal{F}}(T^*N)$. As a consequence of Lemma \ref{inj} we have $\left\{f,\mathcal{J}_{\mathcal{F}}\right\}_{T^*M}\subset\mathcal{J}_{\mathcal{F}}$, which together with Lemma \ref{normif} gives $f\in N(\mathcal{I}_\mathcal{F})$.
\end{proof}

\appendix

\section{Almost Killing Lie algebroids}\label{app}
In the appendix we recall the notion of almost Killing Lie algebroids as defined previously in  \cite{stroblkotov} and provide their relation to module SRFs defined in this paper. (See, in particular, Proposition \ref{prop:almostKilling} below, but also Theorem
\ref{thm:almostKilling} in the main text).

\begin{Definition} A vector bundle $A\to M$ equipped with a  vector bundle morphism $\rho\colon A\to TM$ covering the identity is called an \emph{anchored bundle}. Let $(A,\rho)$ be an anchored bundle  equipped with a skew-symmetric bracket $[\cdot,\cdot]_A$ on $\Gamma(A)$. The triple $(A,\rho,[\cdot,\cdot]_A)$ is called an \emph{almost Lie algebroid} if the induced map $\rho\colon \Gamma(A)\to \mathfrak{X}(M)$ preserves the brackets, and the Leibniz rule is satisfied:

\begin{equation*}
    [s,fs']_A=(\rho(s)\cdot f)\,s'+f[s,s']_A.
\end{equation*}
\end{Definition}

\begin{Definition}
Let $(A,\rho)$ be an anchored bundle over $M$ and $E\to M$ a vector bundle over the same base. An \emph{$A$-connection} on $E$ is a $C^{\infty}(M)$-linear map $^A\nabla$ from $\Gamma(A)$ to $\mathrm{Hom}_\R(\Gamma(E),\Gamma(E)))$ satisfying

\begin{equation*}
    ^A\nabla_s(fe)=(\rho(s)\cdot f)\, e+f^A\nabla_se,
\end{equation*}
for every $f\in C^{\infty}(M)$, $e\in \Gamma(E)$ and $s\in \Gamma(A)$.
\end{Definition}

\noindent An anchored bundle $(A,\rho)$ together with an ordinary connection  on $A$, $\nabla\colon \Gamma(A)\to \Gamma(T^*M\otimes A)$, defines an A-connection $^A\nabla$ on $TM$ by:

\begin{align}
    ^A\nabla_sX&:=\mathcal{L}_{\rho(s)}X+\rho(\nabla_Xs),\label{tau}
\end{align}
valid for every $s\in \Gamma(A)$ and $X\in \mathfrak{X}(M)$. Note that by assuming the Leibniz rule and the commutativity of $^A\nabla_s$ with contractions, these derivations can be extended to arbitrary tensor powers of  $TM$ and $T^*M$.

\begin{Definition}
Let $(A,\rho,[\cdot,\cdot]_A)$ be an almost Lie algebroid over a Riemannian manifold $(M,g)$ and $\nabla\colon \Gamma(A)\to \Gamma(T^*M\otimes A)$ a connection on $A$. Then $(A,\nabla)$ and $(M,g)$ are called compatible if

\begin{equation*}
    ^A\nabla g=0,
\end{equation*}
where the $A$-connection $^A\nabla$ is defined by Equation \eqref{tau}. The triple $(A,\nabla,g)$ is called a $\emph{Killing almost Lie algebroid}$ over $M$.
\end{Definition}

\begin{Lemma}\label{lietau}
Let $(A,\rho)$ be an anchored vector bundle over a Riemannian manifold $(M,g)$, and let $\nabla$ be an ordinary connection on $A$. The triple $(A,\nabla,g)$ satisfies $^A\nabla g=0$ if and only if for every $X,Y\in \mathfrak{X}(M)$ and $s\in \Gamma(A)$ we have

\begin{equation*}
    (\mathcal{L}_{\rho(s)}g)(X,Y)=g(\rho(\nabla_Xs),Y)+g(X,\rho(\nabla_Ys)).
\end{equation*}
\end{Lemma}

\noindent \begin{proof}
By Equation \eqref{tau}, for every vector field $X\in\mathfrak{X}(M)$
\begin{align*}
    (^A\nabla_s g)(X,X)&=^A\nabla_s(g(X,X))-2g(^A\nabla_sX,X)\\
                            &=(\mathcal{L}_{\rho(s)}g)(X,X)-2g(\rho(\nabla_Xs),X).
\end{align*}
Consequently, $^A\nabla g=0$ if and only if

\begin{equation*}
    (\mathcal{L}_{\rho(s)}g)(X,X)=2g(\rho(\nabla_Xs),X) . 
\end{equation*}
\end{proof}

\begin{Proposition}\label{prop:almostKilling}
Let $(M,\mathcal{F})$ be an SF on a Riemannian manifold $(M,g)$. Then the triple $(M,g,\mathcal{F})$ is a module SRF if and only if it is locally generated by Killing almost Lie algebroids, i.e. $\forall q\in M$, there exist an open neighborhood $U\in M$ containing $q$ and a Killing almost Lie algebroid $(A_U,\nabla,g_U)$ over $(U,g_U)$ such that $\rho\left(\Gamma_c(A_U)\right)=\iota_U^{-1}\mathcal{F}$.
\end{Proposition}

\noindent \begin{proof}
Assume that $(M,g,\mathcal{F})$ is a module SRF and $q\in M$. Then there exists an open neighborhood $U\in M$ containing $q$ such that $\iota_U^{-1}\mathcal{F}$ is generated by finitely many vector fields $V_1,\ldots,V_N\in \mathfrak{X}(U)$ for some positive integer $U$. By involutivity of $\iota_U^{-1}\mathcal{F}$, the trivial vector bundle $A_U$ of rank $N$ with a frame $e_1,\ldots,e_N\in \Gamma(A_U)$ together with the anchor map $\rho\colon A_U\to TM$, $e_a\mapsto V_a$ for $a=1,\ldots,N$, can be equipped with an almost Lie algebroid structure. By Lemma \ref{srfgen} there exist $1$-forms $\omega_a^b\in \Omega^1(U)$ such that 

\begin{equation*}
    \mathcal{L}_{V_a}g=\sum_{b=1}^N\omega_a^b\odot \iota_{V_b}g\,\,\,\,\,\forall a,b=1,\ldots,N .
\end{equation*}
Now if we define $\nabla e_a=\sum_{b=1}^N\omega_a^b\otimes e_b$, for every $X,Y\in \mathfrak{X}(U)$, we have

\begin{align*}
    (\mathcal{L}_{\rho(e_a)}g)(X,Y)&=(\mathcal{L}_{V_a}g)(X,Y)\\
    &=\sum_{b=1}^N\left((\iota_X\omega_a^b)g(V_b,Y)+(\iota_Y\omega_a^b)g(X,V_b)\right)\\
    &=g\left(\rho\left(\sum_{b=1}^N\iota_X\omega_a^be_b\right),Y\right)+g\left(X,\rho\left(\sum_{b=1}^N\iota_Y\omega_a^be_b\right)\right)\\
    &=g(\rho(\nabla_Xe_a),Y)+g(X,\rho(\nabla_Ye_a)).
\end{align*}
Consequently, by Lemma \ref{lietau}, $(A_U,\nabla_U,g_U)$ is a Killing almost Lie algebroid and we have $\rho(\Gamma_c(A_U))=\iota_U^{-1}\mathcal{F}$.
Conversely, Assume that $(M,\mathcal{F})$ is locally generated by Killing almost Lie algebroids. Let $q\in M$, and take a neighborhood $U\in M$ containing $q$ with a Killing almost Lie algebroid $(A_U,\nabla,g_U)$ over $(U,g_U)$ such that $\rho(\Gamma_c(A_U))=\iota_U^{-1}\mathcal{F}$. By choosing $U$ small enough, we can assume that $A_U$ is trivial and there is a global frame $e_1,...,e_N\in \Gamma(A_U)$. Then there exist 1-forms $\omega_a^b\in \Omega^1(U)$ such that 

\begin{equation*}
    \nabla e_a=\sum_{b=1}^N\omega_a^b\otimes e_b\,\,\,\,\,\forall a,b=1,\ldots,N.
\end{equation*}
With $V_a:=\rho(e_a)$ for $a,b=1,\ldots,N$, by Lemma \ref{lietau}, for every $X,Y\in\mathfrak{X}(U)$ one has 

\begin{align*}
    (\mathcal{L}_{V_a}g)(X,Y)&=\sum_{b=1}^N\left((\iota_X\omega_a^b)g(V_b,Y)+(\iota_Y\omega_a^b)g(X,V_b)\right)\\
    &= \left(\sum_{b=1}^N\omega_a^b\odot \iota_{V_b}g\right)(X,Y).
\end{align*}
This implies, using Lemma \ref{srfgen} and Proposition \ref{locality}, that $(M,g,\mathcal{F})$ is a module SRF.
\end{proof}


\begin{ack}
We are grateful to Anton Alekseev, Camille Laurent-Gengoux, Ricardo Mendes, Leonid Ryvkin  and, in particular, to Marco Zambon, Mateus de Melo and Serge Parmentier 
for stimulating discussions related in one way or another to the present subject.
\end{ack}

\begin{funding}
This work was supported by the LABEX MILYON (ANR-10-LABX-0070) of Universit\'e de Lyon, within the program ``Investissements d'Avenir'' (ANR-11-IDEX-0007) operated by the French National Research Agency (ANR). We also acknowledge having profited from the marvellous environment provided  within the program ``Higher structures and Field Theory'' at the ERWIN SCHR\"ODINGER INSTITUTE in Vienna.
\end{funding}


\end{document}